\documentclass[12pt]{amsart}
\usepackage{amscd,amssymb,times,fullpage,graphicx}
\newtheorem{thm}{Theorem}
\newtheorem{prop}[thm]{Proposition}
\newtheorem{lem}[thm]{Lemma}
\newtheorem{cor}[thm]{Corollary}

\newtheorem{ex}[thm]{Example}

\theoremstyle{definition}
\newtheorem{defn}[thm]{Definition}

\DeclareMathOperator{\SO}{SO}

\DeclareMathOperator{\coh}{coh}

\DeclareMathOperator{\SU}{SU}

\DeclareMathOperator{\KO}{KO}

\newcommand{\C}{\mathbb{ C}}

\newcommand{\Z}{\mathbb{ Z}}

\newcommand{\R}{\mathbb{ R}}

\newcommand{\im}{\rm im}

\begin{document}

\title{Positive scalar curvature with symmetry}

\author{Bernhard Hanke}

\begin{abstract} We show an equivariant bordism principle for constructing metrics of positive scalar curvature 
that are invariant under a given group action. Furthermore, we develop a new codimension-$2$ surgery technique 
which removes singular 
strata from fixed point free $S^1$-manifolds while preserving equivariant 
positive scalar curvature. These results are applied to derive
the following theorem: 
Each closed fixed point free $S^1$-manifold
of dimension at least $6$ whose isotropy groups have odd order and whose union of maximal
 orbits is simply connected and 
not spin, carries an $S^1$-invariant metric of positive scalar curvature.   
\end{abstract}

\date{\today; MSC 2000: primary 53C23, 53C25, 57R85; secondary  57R65}

\maketitle

\tableofcontents

\section{Introduction} \label{s:intro}

A classical theme in differential geometry is the investigation of   
topological conditions that are necessary or sufficient for 
the existence of a particular kind of geometric structure 
on a given smooth manifold. In the context of Riemannian 
metrics of positive scalar curvature, this question has revealed a strong 
link between subtle differential topological invariants of  smooth 
manifolds and their geometry. A prominent role in this context is 
played by an effective method for constructing metrics of positive 
scalar curvature described in the seminal work by Gromov-Lawson \cite{GL}
and Schoen-Yau \cite{SY}: The class of 
smooth manifolds admitting metrics of positive scalar curvature is 
closed under surgery in codimension 
at least three. Based on this principle, the existence question 
for positive scalar curvature metrics can be translated into 
a bordism problem that is then discussed with the help of 
powerful algebraic-topological means. The effectiveness of this approach is 
illustrated by the following result, which provides
a complete classification of simply connected closed manifolds  of dimension 
at least $5$ that admit metrics of positive scalar curvature.

\medskip

\noindent {\bf Theorem.} \begin{it} Let $M$ be a closed simply 
connected manifold of dimension at least $5$. 
\begin{itemize}
   \item[i.)] {\rm (Gromov-Lawson \cite{GL})} \label{GL} If $M$ does not admit a spin
structure, then $M$ carries a metric of positive scalar curvature.
   \item[ii.)] {\rm (Lichnerowicz \cite{Li}, Hitchin \cite{Hi}, 
               Stolz \cite{St})} \label{St} If $M$ admits a spin 
             structure, then $M$ carries a metric of positive scalar 
             curvature if and only if $\alpha(M)=0$.
\end{itemize}
\end{it}

\medskip

Here $\alpha(M)\in \KO_n$ is an invariant (defined in \cite{Hi}) 
closely related to $\hat{A}(M)$. For non-simply connected manifolds,
there are refined versions of this 
obstruction with values in the $K$-theory of certain $C^*$-algebras
associated to the fundamental group of $M$. A complete classification 
of manifolds admitting 
metrics of positive scalar curvature is not known, even in the 
case of finite fundamental groups.   

We will investigate the positive scalar curvature 
question in an equivariant context: Given a closed smooth manifold
$M$ equipped with a smooth action  
of a compact Lie group $G$, does $M$ 
admit a positive scalar curvature metric which is 
invariant under the $G$-action? Even if $M$ admits 
a (nonequivariant) metric of positive scalar curvature, this may be 
a nontrivial problem:  B\'erard Bergery shows in \cite{Bebe}, Example 9.1., 
that 
averaging a positive scalar 
curvature metric over a group of symmetries may destroy the 
positivity of the scalar curvature.

If $G$ is finite and the action on $M$ is free, this problem is   
equivalent to asking whether $M / G$, a closed manifold 
whose fundamental group is a certain extension of $G$ (assuming 
that $M$ is connected), admits a metric of positive scalar 
curvature -   and this is a nonequivariant question. Another extreme case is 
that of (not necessarily free) actions of compact Lie groups whose
identity components are non-abelian. This problem is
completely settled by the following result. 

\medskip

\noindent {\bf Theorem.} (Lawson-Yau \cite{LawYau}) \begin{it} 
If a compact connected manifold $M$ is equipped with 
an effective action of a compact Lie group $G$ whose 
identity component is non-abelian, then it admits a Riemannian metric 
of  positive scalar curvature which is invariant under the given $G$-action. 
\end{it}

\medskip

In its original form, this theorem only states that 
$M$ carries a Riemannian metric of positive scalar curvature, 
but one easily checks that the construction in {\it loc.~cit.}~yields 
a metric which is in fact also $G$-invariant. 

The proofs of the Lawson-Yau 
theorem and of the Gromov-Lawson and Stolz   
theorems differ in an 
essential way: The first one yields  an explicit 
metric of positive scalar curvature. In contrast,
the proofs of the latter are based on structure results for
the oriented and spin bordism rings which rely 
on homotopy theoretic considerations. In particular, they 
do not provide a direct description of the positive scalar 
curvature metrics in question.

Important aspects of the equivariant positive 
scalar curvature problem 
were discussed by B\'erard Bergery in \cite{Bebe}. 
In this paper the Kazdan-Warner trichotomy, the 
Yamabe problem, the index obstruction to positive scalar curvature 
and the surgery principle for constructing metrics
of positive scalar curvature are formulated in an equivariant context.

Further elaborations of the index theoretic obstruction in 
the case of $S^1$-manifolds are carried out by Lott
in \cite{Lott}. Rosenberg-Weinberger  \cite{RosWei} provide an  
ad-hoc discussion of  an equivariant bordism 
principle for constructing metrics of positive scalar curvature 
on simply connected manifolds which are invariant under spin preserving $\Z/p$-actions, 
see  Theorem 2.3. in {\em loc.~cit}. However, 
the proof contained in this paper requires more assumptions
than stated in the theorem (see the discussion following 
Corollary \ref{zetpe} below). Relying  
on this (potentially problematic) result, Farsi \cite{Far} has  discussed some 
instances of the equivariant positive scalar curvature problem 
on spin manifolds of low dimension equipped with $\Z/p$-actions.

Our work is devoted to a 
systematic exploration 
of the surgery and bordism techniques for  
constructing equivariant positive scalar
curvature metrics. Part of 
the nonequivariant discussion can be translated 
more or less directly to the equivariant context. This applies 
in particular to the surgery principle of 
Gromov-Lawson and Schoen-Yau.  The paper \cite{Bebe} 
formulates an equivariant analogue of this 
fundamental result (see \cite{Bebe}, Theorem 11.1) without proof. 
In Section \ref{equivGL} of our work, we will recapitulate the 
essential steps of the argument  in \cite{GL} 
and will explain how they translate to an equivariant setting. 

The following Section \ref{bordism} is devoted to a proof 
of the first main result, 
a general bordism principle for constructing equivariant positive
scalar curvature metrics. 

\medskip

\noindent {\bf Theorem A}. \begin{it} Let $Z$ be a compact connected 
oriented $G$-bordism between the closed $G$-manifolds $X$ and $Y$. Assume the following: 
\begin{itemize}
       \item[i.)] The cohomogeneity of $Z$ is at least $6$,  
       \item[ii.)] the inclusion of maximal orbits 
                $Y_{max} \hookrightarrow Z_{max}$ is a 
                nonequivariant $2$-equivalence (i.e.~a bijection 
                 on $\pi_0$, an isomorphism on $\pi_1$ and a 
                surjection on $\pi_2$),
        \item[iii.)] each singular stratum of codimension $2$ in $Z$ meets $Y$. 
\end{itemize}
Then, if $X$ admits a $G$-invariant metric of positive scalar
curvature, the same is true for $Y$. 
\end{it}

\medskip

This statement makes clear that the bordism 
principle for constructing equivariant metrics of positive 
scalar curvature metrics does definitely not require a strong gap hypothesis 
as envisaged in \cite{RosWei}, Remark 2.4, where the authors assume that 
for any two closed subgroups $H,K \subset G$ with $K \subset H$, the codimension 
of each component of $Z_{(H)}$ contained in the closure of $Z_{(K)}$ is  
either $0$ or at least $3$ in $Z_{(K)}$ (here $Z_{(H)}$ denotes 
the set of points in $Z$ whose 
isotropy groups are conjugate to $H$).

Our Theorem A is almost a direct analogue of the corresponding nonequivariant 
result (see \cite{StolzCongress}, Theorem 3.3). In particular 
the dimension restriction i.) and 
the connectivity restriction for the inclusion $Y_{max} 
\hookrightarrow Z_{max}$ stated in point ii.) translate to 
analogous requirements in the nonequivariant setting if  $G=\{1\}$. 
However, if $G$ is not trivial, we need an additional assumption on codimension-$2$ singular 
strata. The plausibility of such an assumption is illustrated in  Proposition 
\ref{illustrate}: Any  equivariant handle decomposition of 
$\Z/p$-bordisms $Z$ with codimension-$2$ fixed point components 
disjoint from $Y$ contain handles of codimension $0$ or $2$
(independent of the connectivity of the map $Y_{max} \hookrightarrow Z_{max}$). 
This points towards a fundamental limitation of the method of 
equivariant handle decompositions to construct 
equivariant metrics of positive scalar curvature.   

Theorem A is useful for constructing equivariant metrics of 
positive scalar curvature only if it can be combined
with powerful structure results for geometric equivariant bordism 
groups, which  imply that the manifold $X$ in Theorem A can 
be assumed to admit an equivariant positive scalar 
curvature metric under some general assumptions on the
manifold $Y$. Two main difficulties occur at this point. 
Firstly, explicit geometric generators 
of equivariant bordism groups are known only 
in a very limited number of cases.  Secondly, 
whereas conditions i.) and ii.)  in Theorem A 
can be achieved under fairly general assumptions
on the manifold $Y$ (by performing appropriate surgeries on 
$Z_{max}$ - cf. the second proof of Proposition \ref{free}), it is a 
priori not clear under 
what circumstances condition iii.) holds. 

In sections \ref{genial} and \ref{anwend}, 
we shall present a way to avoid the difficulties inherent in 
condition iii.) if $G =S^1$ and the 
$G$-action on $Z$ is fixed point free. The idea we use  
is  to alter a given bordism $Z$ by cutting out equivariant 
tubes connecting $Y$ with each of the codimension-$2$ 
singular strata in $Z$ that are disjoint from $Y$. This replaces 
the bordism $Z$ and the manifold $Y$ by other  manifolds $Z'$ 
and $Y'$ so that each codimension-$2$ singular stratum in $Z'$ 
meets $Y'$. In particular, Theorem A can be applied to $Z'$ (after 
some more manipulations of $Z'$, but we omit these details here). We must now 
understand 
how $Y$ can be recovered from $Y'$. 
A closer inspection of the situation shows that $Y'$ is obtained from $Y$ by adding 
certain codimension-$2$
singular strata with finite isotropies. Conversely, $Y$ can be reconstructed from 
$Y'$ by a Dehn-like codimension-$2$ surgery process that removes 
these additional singular strata and puts back free ones instead. The examination 
of this surgery step is the content of Section \ref{genial}. In Theorem 
\ref{heikel} of this section we show by a somewhat involved geometric 
argument that this surgery step preserves 
the existence of $S^1$-invariant positive scalar curvature metrics 
under fairly general assumptions.  
Roughly speaking, we replace the ``bending outwards'' process in the surgery 
step due to Gromov-Lawson and 
Schoen-Yau by a ``bending inwards'' process. We remark 
that this kind of positive scalar curvature preserving codimension-$2$ surgery  only 
works  under the additional $S^1$-symmetry. 

Arguing in this rather roundabout manner, assumption iii.) of Theorem A 
is no longer a true obstacle against the  
construction of equivariant positive 
scalar curvature metrics on fixed point free $S^1$-manifolds.
In combination with a classical theorem of Ossa \cite{Ossa}, 
which states that fixed point free oriented $S^1$-manifolds satisfying 
condition $\C$ (cf. Definition \ref{conditionce}) are oriented 
$S^1$-boundaries, we get the following equivariant version 
of the Gromov-Lawson theorem stated above.    

\medskip

\noindent {\bf Theorem B.} \begin{it}  Let $M$ be a closed
fixed point free $S^1$-manifold 
satisfying condition $\C$  and of cohomogeneity at 
least $5$. 
If the union of maximal orbits of $M$  
is simply connected and  does not 
admit a spin structure,  then $M$ admits an $S^1$-invariant metric of positive scalar curvature.
\end{it}

\medskip

By Lemma \ref{oddorder}, the manifold $M$ satisfies condition $\C$, if all isotropy groups 
have odd order. We remark that no additional assumption on codimension-$2$ singular strata in $M$ 
is needed in Theorem B. It is not clear at present to what extent Ossa's theorem 
can be generalized to the spin case, so that we will not 
discuss the corresponding $S^1$-equivariant analogue of Stolz' theorem in this paper.

\medskip

This paper forms part of my habilitation thesis at the 
Ludwig-Maximilians-Universit\"at M\"unchen. 
I am grateful to D. Kotschick for constant help and encouragement, to 
L. B\'erard Bergery for sending me a copy of \cite{Bebe} and 
to J. Rosenberg, T. Schick and S. Stolz for useful comments.  
The work on this paper was supported by a 
research grant within the DFG emphasis program ``Global Differential Geometry''.

\section{Equivariant surgery theorem} \label{equivGL}

\begin{defn} Let $G$ be a compact Lie group. An 
{\em equivariant $G$-handle} is  a $G$-space of the form 
\[
     G \times_H ( D(V) \times D(W) ) 
\]
where $H \subset G$ is a closed subgroup and 
$D(V)$ and $D(W)$ are unit discs in some orthogonal $H$-representations
$V$ and $W$. 
We call $\dim V$ the {\em dimension} and $\dim W$ the 
{\em codimension} of the $G$-handle in question. 
\end{defn}

If $Z$ is a smooth $G$-manifold with boundary and 
\[
   \phi :  G \times_H (S(V) \times D(W))  \hookrightarrow \partial Z 
\]
is a $G$-equivariant embedding, we can glue the given 
$G$-handle  to $Z$ along $\phi$ and obtain a $G$-space 
\[
    Z' =  Z \cup_{\phi} \big( G \times_H (D(V)  \times D(W)) \big)
\]
which can canonically be equipped with the structure of a smooth $G$-manifold 
by straightening corners. We say that $Z'$ is obtained
from $Z$ by {\em attaching} the given $G$-handle. Correspondingly 
we call the subset 
\[
      \phi (G \times_H (S(V) \times 0)) \subset \partial Z
\]
the {\em attaching sphere} of the $G$-handle. If $M$ 
is a smooth $G$-manifold, $Z := M \times (0,1]$ and 
$Z'$ is obtained from $Z$ by attaching a $G$-handle of dimension 
$d\geq 0$ and codimension $c\geq 0$ to $\partial Z$ (which can be identified 
with $M$), we say that  $\partial Z'$ is obtained from $M$ by 
{\em equivariant surgery} of {\em dimension} 
$d-1$ and {\em codimension} $c$ (in the case when $G$ is a finite 
group, these numbers refer to the dimension and 
codimension of the attaching sphere in $M$). In this terminology, 
equivariant surgery of dimension $-1$ (i.e. $d=0$) 
amounts to adding a disjoint component of the form 
\[
    G \times_H (D^0 \times S(W)) 
\]
to $M$. 

The equivariant generalization of the surgery result due to Gromov-Lawson
and Schoen-Yau is fairly straightforward as was already remarked in
\cite{Bebe}.

\begin{thm}[cf.~\cite{Bebe}, Theorem 11.1] \label{equivgl}
Let $M$ be a 
(not necessarily compact) $G$-manifold equipped with a 
$G$-invariant 
metric of positive scalar curvature and  let $N$ 
be obtained from $M$ by equivariant surgery of 
codimension at least $3$. Then $N$ carries a $G$-invariant 
metric of positive scalar curvature which 
can be assumed to coincide with the given metric on $M$ outside 
a prescribed  neighbourhood (which may be arbitrarily small)
around the attaching sphere in $M$.  
\end{thm}

\begin{proof} We start with an equivariant embedding 
\[
   \phi :  G \times_H (  S(V) \times D(W))  \hookrightarrow M \, . 
\]
After blowing up the metric on $M$ if we wish to increase
the injectivity radius of $M$, we can (up to $G$-isotopy) assume that the radial lines 
in $D(W)$ are mapped to unit speed geodesics in $M$ that are orthogonal 
to the attaching sphere 
\[
    \phi(G \times_H( S(V) \times 0)) \subset M \, . 
\]
The generalization of the construction 
in \cite{GL} to the  equivariant situation is possible
because the crucial part of the argument
uses the distance from the attaching sphere in $M$ and 
in the equivariant context,
this is automatically a $G$-equivariant  function. 

For completeness and because \cite{Bebe} contains 
only a sketch of the proof, we summarize the essential 
steps of the construction.

If the $G$-handle $G \times_H (D(V) \times D(W))$ is of dimension 
$0$ (this case - which may very well occur - is usually skipped in 
the literature), then the corresponding surgery step amounts 
to adding a new component of the form 
\[
    G \times_H (D^0 \times S(W)) 
\]
to $M$. This new component  can be written as the total space of a 
$G$-equivariant fibre bundle
\[
    S(W) \hookrightarrow G  \times_H S(W) \to G/H 
\]
and an application of the O'Neill formula together with the assumption 
that $\dim S(W) \geq 2$ shows that this new component admits 
a $G$-invariant metric of positive scalar curvature. 

From now on, we assume that $\dim V \geq 1$. Let $g$ be the 
given metric on $M$.    
We consider the pull back metric $\phi^*(g)$ on 
\[
    G \times_H ( S(V)  \times D(W)) \, .  
\]
For $0 < \epsilon < 1$ let $\phi^*(g)_{\epsilon}$
denote this pull back metric restricted to 
\[
    G \times_H (S(V) \times S_{\epsilon}(W))
\]
where $S_{\epsilon}(W)$ is the $\epsilon$-sphere in $W$ (with 
respect to the given euclidean metric on $W$). Note that we have a $G$-fibre bundle
\[
    S_{\epsilon}(W) \hookrightarrow G \times_H (S(V) \times S_{\epsilon}(W)) \to G \times_H (S(V) \times 0)  
\]
which is nothing but the $\epsilon$-sphere bundle (with respect 
to the  metric $\phi^*(g)$) of the fibration 
\[
 D(W) \hookrightarrow G \times_H (S(V) \times D(W)) \to G \times_H S(V) \, .  
\]
However, the metric 
$\phi^*(g)_{\epsilon}$ on $G \times_H (S(V) \times S_{\epsilon}(W))$ will not, 
in general, be a Riemannian submersion metric. 

We now construct a Riemannian submersion metric $h_{\epsilon}$ on 
this bundle which is 
uniquely characterized by the following properties:  
\begin{itemize}
   \item[i.)] On the base $G \times_H S(V)$ it is the restriction of $\phi^*(g)$, 
   \item[ii.)] on the fibres it is the usual round metric on $S_{\epsilon}(W)$, 
   \item[iii.)] the horizontal subspaces are determined by the normal connection 
         of the embedding $\phi\big( G \times_H (S(V) \times 0) \big) \subset 
         M$. 
\end{itemize}

The fibres 
$S_{\epsilon}(W)$ are totally geodesic with respect to $h_{\epsilon}$.
We pick a $G$-invariant metric on $G \times_H  D(V)$  
which near the boundary is the product metric 
$\phi^*(g)|_{G \times_H S(V)} \oplus  dr^2$ (where $r$ denotes the radial coordinate in $D(V)$). Using this, 
we construct a $G$-invariant submersion metric $\overline{h}_{\epsilon}$ 
on the bundle
\[
    S_{\epsilon}(W) \hookrightarrow G \times_H (D(V) \times S_{\epsilon}(W)) 
   \to G \times_H (D(V) \times 0)  
\]
which is a $G$-equivariant product $h_{\epsilon} \oplus  dt^2$ near the boundary 
and has totally geodesic fibres $S_{\epsilon}(W)$. Because $\dim S(W) \geq 2$,
an application of the O'Neill formula shows that 
there is a small $\epsilon_0 > 0$ such that $\overline{h}_{\epsilon}$
has positive scalar curvature for all $0 < \epsilon < \epsilon_{0}$.

The following fact is the equivariant version of  the crucial ``bending outwards'' 
process described in the work of Gromov-Lawson \cite{GL} (see also the elaborations of 
this argument in \cite{RS} and \cite{Wehrheim}). Again, the proof in the 
equivariant case does not lead to any further complications so that 
we confine  ourselves to a clear statement of the relevant fact. 

\medskip

\noindent \begin{it} There is a $T \in \R^+$ and 
a $G$-invariant metric $\gamma$ of positive scalar curvature
on 
\[
     G \times_H (S(V) \times S(W))  \times [0,T]
\]
which on the left hand part $G \times_H (S(V) \times S(W)) \times [0, \delta)$ 
(where $\delta > 0$ is some small number),
is isometric to  $\phi^*(g)$ restricted to a collar $\delta$-neighbourhood
of $G \times_H (S(V) \times S(W))$ in 
$G \times_H(S(V) \times D(W))$ and on the right hand part
$G \times_H (S(V) \times S(W)) \times (T - \delta, T]$, is  
a Riemannian product metric $h_{\epsilon} \oplus dt^2$. 
Here $\epsilon$ can be chosen arbitrarily in some interval $(0, \epsilon_1)$ where $\epsilon_1 > 0$ 
is an appropriately chosen small constant. 
\end{it}

\medskip

We now pick some  $\epsilon$ in the interval $(0, \min(\epsilon_0, \epsilon_1))$ and 
glue the piece
\[
    \big( G \times_H (S(V) \times S(W)) \times [0,T], \gamma \big) 
\]
to 
\[
       \big( M \setminus \phi\big( G \times_H (S(V) \times D(W)) \big), g  \big)
\]
along $G \times_H (S(V) \times S(W)) \times \{0\}$ and then glue the Riemannian manifold 
\[
   \big( G \times_H (D(V) \times S_{\epsilon}(W)), \overline{h}_{\epsilon} \big)
\]
along $G \times_H (S(V) \times S_{\epsilon}(W))$ 
to the boundary of the resulting space. The manifold thus obtained is $G$-diffeomorphic 
to $N$ and carries a $G$-invariant metric of positive scalar curvature by construction. 
\end{proof}

\section{Equivariant bordism theorem} \label{bordism}

The reduction of the construction of  (nonequivariant) positive 
scalar curvature metrics to a bordism problem relies on the 
fact that a bordism $Z$ (which is always assumed to be compact)
between two closed manifolds $X$ and $Y$ 
can be decomposed  into a series of handle attachments and that 
the (co-)dimensions of these handles can be controlled in 
terms of the connectivity of the map $Y  \hookrightarrow Z$.
This uses the existence of Morse functions on $Z$ and the technique 
of handle cancelations. For an exposition of these methods, 
see e.g.~\cite{Kosinski, Lueck, Milnorm, Milnorh}. 
If the inclusion $Y \hookrightarrow Z$ is a $2$-equivalance, 
this implies that $Y$ can be obtained from $X$ by a series 
of surgeries of codimension at least $3$ so that 
$Y$ carries a metric of positive scalar curvature if $X$ carries such 
a metric.

Not all of these steps carry over directly to the 
equivariant case. Indeed, it was observed by Wassermann \cite{Was}  that 
Morse theory can be formulated in an equivariant 
setting and can be used to construct decompositions 
of closed smooth $G$-manifolds into $G$-handles. 
However, it is well known that handle cancelation
does not work in the equivariant context in full generality.
This leads to counterexamples to equivariant 
analogues of the  h- and s-cobordism 
theorems (see e.g.~\cite{Kawa, Kawa2}) and is also a main obstacle 
against translating Wall's surgery theory to an equivariant 
setting. In order to circumvent these difficulties,
one usually works with gap hypotheses of the 
form that each singular stratum which is properly contained in the 
closure of another singular stratum $F$ must have a large enough
codimension in $F$ (cf. \cite{Kawakubo}) or one formulates the $s$-cobordism 
theorem in an isovariant context, see e.g. \cite{Lu}, Theorem 4.42. In some sense, 
Theorem A follows this isovariant viewpoint.  

Our proof of this result is based on two 
observations. The first one is that we need to cancel only $G$-handles in 
$Z$ of codimension less than $3$ so that the full power of a 
handle cancelation machinery is not necessary. The second 
 - Lemma \ref{codestimate} below - is that the 
codimensions of $G$-handles occurring 
in $Z$ are related to the codimensions of the singular strata to  
which they are attached if the handle decomposition 
of $Z$ is induced by a $G$-Morse function which is
special in the sense of \cite{Mayer}. 

We start by recalling some important notions from equivariant differential topology. 
Let $G$ be a compact Lie group and let 
$M$ be a smooth $G$-manifold. For a closed subgroup $H \subset  G$, we denote by 
$(H)$ the conjugacy class of $H$ in $G$. The set of 
conjugacy classes of subgroups of $G$ is partially ordered
by writing $(H) \leq (K)$ if and only if $H$ is conjugate 
to a subgroup of $K$. For  $x \in M$ 
let $G_x \subset G$ be the isotropy group of $x$. Furthermore,  
we use the following notation: 
\begin{itemize}
   \item[i.)] $M_{(H)} : = \{ x\in M ~|~ (G_x) =  (H) \}$, 
   \item[ii.)] $M^H : = \{ x \in M ~|~ H \subset G_x\} = \{ x \in M ~|~ hx = x {\rm 
                     ~for~all~} h \in H \}$. 
 \end{itemize}
The space $M^H$ is a closed submanifold of $M$, but usually consists of 
different orbit types and is in general not $G$-invariant unless $G$ is abelian. 

In contrast, $M_{(H)}$ is an (in general not closed) 
$G$-submanifold of $M$. The space $M_{(H)}$ is called the 
{\em $H$-orbit bundle} of $M$.
It consists of all points in $M$ with isotropy groups conjugate to $H$
and these form exactly those $G$-orbits which are $G$-diffeomorphic 
to the left $G$-space $G /H$.  There is a $G$-fibre bundle  
\[
       G/H \hookrightarrow M_{(H)} \to M_{(H)} / G  \, . 
\]
We cite

\begin{prop}[\cite{tD}, Theorem (5.14)] Suppose $M$ is a 
$G$-manifold and $M/G$ is connected. Then there exists a unique 
isotropy type $(H)$ such that $M_{(H)}$ is open and dense in $M$. 
The space $M_{(H)}/G$ is connected. Each isotropy type 
$(K)$ satisfies $(H) \leq (K)$. The set $M^H$ intersects each orbit. 
\end{prop}

We call the space $G/H$ with $H$ as in the last proposition the 
{\em principal (or maximal) orbit type} of $M$. Accordingly we call 
$(H)$ the {\em minimal isotropy type}. From now on
we denote the minimal isotropy type by $(H_{min})$ and use the 
shorthand notation 
\[
     M_{max} : = M_{(H_{min})}   
\] 
for the union of maximal orbits in $M$. This subset is open 
and dense in $M$. 
If $(H) \neq (H_{min})$, we call each component of $M_{(H)}$ a 
{\em singular stratum} of the $G$-action. The {\em cohomogeneity} 
\[
   {\rm coh}(M,G)  
\]
of a connected $G$-manifold $M$ is the codimension of a principal orbit 
in $M$ (this does not depend on the principal orbit chosen). Note the equality
\[
    {\rm coh}(M,G) = \dim M_{max} / G \, . 
\] 

The following example illustrates how different singular strata 
in a given $G$-manifold can be related to each other.

\begin{ex} \label{typical} Let $V_{2}$ and $V_{3}$ be the 
irreducible one dimensional complex 
$\Z/6$-representations of weights $2$ and $3$ respectively. Then the $S^1$-manifold 
\[
  M := S^1 \times_{\Z/6}  ( V_{2} \times V_{3} ) 
\]
has three  singular strata with isotropy groups $\Z/2$, $\Z/3$ and $\Z/6$. The 
singular strata with isotropy $\Z/2$ and $\Z/3$ are not closed in $M$. More 
precisely, the intersection of the closures of these two singular strata 
is the third singular stratum with 
isotropy $\Z/6$. The maximal orbit type of $M$ is equal to $S^1$ (in other words, the 
given $S^1$-action is effective) and the cohomogeneity of $M$ is 
equal to $4$. 
\end{ex}

Now we give  those notions and results  from equivariant Morse theory
that are important for our discussion. We will mainly 
refer to the papers of Mayer \cite{Mayer} and Wassermann \cite{Was}.

\begin{defn} Let $M$ be a closed $G$-manifold 
and let 
\[
   f:M \to \R
\]
be a smooth $G$-equivariant map where $\R$ is equipped with 
the trivial $G$-action. 
An orbit 
\[
   G / H \approx \mathcal{O}  \subset M
\]
is called {\em critical} if for one (and hence any) point $x \in 
\mathcal{O}$ the differential $D_x f$ is zero. The critical 
orbit $\mathcal{O}$ is called {\em nondegenerate} if for each 
$x \in \mathcal{O}$ the following holds: Let
$N_x \subset M$ be a normal slice of $\mathcal{O}$ at $x$. Then 
the $H$-invariant function
\[
     f|_{N_x} : N_x \to \R  
\]
has a nondegenerate critical point at $0$, i.e.~the Hessian 
of $f|_{N_x}$ is nondegenerate. The index (resp.~coindex) of the Hessian 
of $f|_{N_x}$ at $0$ is called the {\em index} (resp.~{\em coindex})
of $f$ at the 
nondegenerate critical orbit $\mathcal{O}$. Note that 
the property of $f$ being 
nondegenerate at $\mathcal{O}$ and the index do not depend on the 
choice of $x \in \mathcal{O}$ or the choice of normal slices. The function $f$ is 
called a {\em $G$-Morse function} if it has only nondegenerate 
critical orbits. 
\end{defn}

If $M$ contains just one orbit type, then $M/G$ is a 
smooth manifold and it is clear that
$G$-Morse functions $M \to \R$ and ordinary Morse functions $M/G \to \R$ 
are in one-to-one correspondence. 

Similar to the nonequivariant case, we have
 
\begin{lem}[\cite{Was}, Density Lemma 4.8.] \label{dense} Let $M$ be 
a closed $G$-manifold. Then the set of $G$-Morse 
functions is dense (and clearly open) in the set of smooth 
$G$-equivariant maps 
\[
    M \to \R
\]
equipped with the $C^{\infty}$-topology. 
\end{lem} 

Passage through critical orbits is described by attaching equivariant
handles: 

\begin{lem}[\cite{Was}, Theorem 4.6.] \label{handledec} Let $M$ be a 
closed $G$-manifold, let 
\[
   f: M \to \R
\]
be a $G$-Morse function and let $[a,b] \subset \R$ contain
exactly one critical value of $f$, lying in $(a,b)$. 
Then $f|_{f^{-1}[a,b]}$ has
finitely many critical orbits 
\[
     \mathcal{O}_1, \ldots, \mathcal{O}_r \, . 
\]
Let $d_1, \ldots, d_r$ be the indices and 
$c_1, \ldots, c_r$ be the coindices of $f$ at these critical 
orbits. Then $f^{-1}((-\infty, b])$ is 
equivariantly diffeomorphic to $f^{-1}((- \infty, a])$ 
with finitely many disjoint $G$-handles 
\[
    G \times_{H_1} (D(V_1) \times D(W_1) ) \, , \ldots \, ,  
    G \times_{H_r} (D(V_r) \times D(W_r))
\]
attached. Here $(H_i)$ is the isotropy type of the orbit 
$\mathcal{O}_i$ and $\dim V_i = d_i$. Furthermore, after
a choice of appropriate coordinates, the Morse 
function $f$ restricted to a slice $D(V_i) \times D(W_i)$ has 
the standard form 
\[
   f(x_1, \ldots, x_{d_i}, y_1, \ldots, y_{c_i}) = f(0) - 
  x_1^2 - \ldots - x_{d_i}^2 + y_1^2 + \ldots + y_{c_i}^2 \, . 
\]
\end{lem}

As an immediate corollary of these results, each closed
$G$-manifold admits a $G$-handle decomposition.

Now let $X$ and $Y$ be  closed $G$-manifolds and let $Z$ be a $G$-bordism 
from $X$ to $Y$, i.e.~$Z$ is a compact $G$-manifold 
whose boundary  splits as 
the disjoint union of $X$ and $Y$. 
In this special case we call a smooth $G$-equivariant map  
\[
   f: Z \to \R
\]
a {\em $G$-Morse function} if in addition to the previous 
requirements it satisfies
\begin{itemize}
   \item[i.)] $f(Z) \subset [0,1]$,  ~ $f|_{X} = 0 $, ~ $f|_{Y} = 1 $, 
   \item[ii.)] the critical values of $f$ are different from $0$ and $1$.
\end{itemize}
For smooth $G$-equivariant maps $Z \to [0,1]$ with  these additional 
properties, analogues of the density Lemma \ref{dense} as well as of 
the passage-through-critical-orbits Lemma \ref{handledec} hold. 

\begin{cor} Let $Z$ be a compact $G$-bordism from $X$ to $Y$ where $X$ and 
$Y$ are closed $G$-manifolds. 
Then the manifold $Z$ can be obtained from $X\times [0,1]$ 
by successively attaching finitely many $G$-handles. 
\end{cor}

If we wish to use this result for constructing $G$-invariant 
positive scalar curvature metrics, we need to  control 
the codimensions of the handles occuring in such a  $G$-handle decomposition. 
This is possible by use of  {\em special} $G$-Morse functions
that were introduced by Mayer in \cite{Mayer}. Note  
that if $f : M \to \R$ is a $G$-Morse function, then 
the restriction of $f$ to every orbit bundle $M_{(H)}$
has the same critical orbits and is also a 
$G$-Morse function (see Satz 3.1 in \cite{Mayer}). 
However, it can be shown by easy examples that 
in general the indices of $f$ and 
of the restricted Morse function do not
coincide. This motivates the following definition. 

\begin{defn}[\cite{Mayer}, Definition 2.1] Let $M$ be a closed $G$-manifold. A $G$-Morse 
function 
\[
   f: M \to \R 
\]
is called {\em special} if for 
each critical orbit $\mathcal{O}$ 
the index of $f$ at $\mathcal{O}$ is equal to the 
index of the restricted $G$-Morse function
\[
    f|_{M_{(H)}} : M_{(H)} \to \R
\]
at $\mathcal{O}$. Here  $(H)$ is the isotropy type of $\mathcal{O}$. 
\end{defn}

For special $G$-Morse functions we have the following genericity statement.

\begin{lem}[\cite{Mayer}, Satz 2.2 and the following Bemerkung] \label{specialdense}
Let $M$ be a closed $G$-manifold. Then the set of special $G$-Morse functions
is dense in the set of all smooth $G$-equivariant functions $M \to \R$ 
in the $C^1$-topology (but in general not in the $C^{2}$-topology). 
\end{lem} 

A similar statement holds for special $G$-Morse functions defined 
on bordisms between two closed $G$-manifolds $X$ and $Y$. 

We get the following variant of  Lemma \ref{handledec}.

\begin{lem} \label{fuenf} Let $M$ be a  closed $G$-manifold, let 
\[
   f: M \to \R
\]
be a special $G$-Morse function and let $[a,b] \subset \R$ contain
exactly one critical value  as before. 
Then $f|_{f^{-1}[a,b]}$ has finitely many  critical orbits 
\[
     \mathcal{O}_1, \ldots, \mathcal{O}_r \, . 
\]
Let $d_1, \ldots, d_r$ be the respective indices of $f$ at these critical 
orbits. Then $f^{-1}((-\infty, b])$ is 
equivariantly diffeomorphic to $f^{-1}((- \infty, a])$ 
with a finite number of disjoint $G$-handles 
\[
    G \times_{H_1} (D^{d_1} \times D(W_1))  \,  ,  \,  \ldots \, , \,  
    G \times_{H_r} (D^{d_r} \times D(W_r))
\]
attached. Here $(H_i)$ is the isotropy type of the orbit 
$\mathcal{O}_i$ and $H_i$ acts trivially on the unit discs $D^{d_i}$ for 
$1 \leq i \leq r$.  
\end{lem}

A similar statement holds for special $G$-Morse functions defined on compact 
$G$-bordisms $Z$. From now 
on we will speak of  $G$-handles of the form $G \times_H (D^d \times D(W))$ 
as {\em special} $G$-handles.

One pleasant feature of  special $G$-Morse functions is that they 
lead directly to $G$-CW structures on the given manifolds, see \cite{Mayer}, Satz 3.3. 

Before we state and prove the main result of this section, we need three
more preparatory lemmas.  The first one says that we have some control on the order
in which $G$-handles occur if we work with special $G$-Morse functions.

\begin{lem} \label{inducedhandledec} Let $Z$ be a compact $G$-bordism as 
before, let 
\[
       Z \to [0,1]
\]
be a special $G$-Morse function and let $H \subset G$ be a closed subgroup. Then 
there is a special $G$-Morse function 
\[
    f : Z \to [0,1]
\]
with the same critical orbits and the same indices on these 
critical orbits as the given $G$-Morse function, but with 
the following additional property: There are two 
noncritical values 
\[
    0 < c < d < 1 
\]
of $f$ so that  for each critical 
orbit $\mathcal{O}$ of $f$ the following equivalences hold:
\begin{eqnarray*}
     f(\mathcal{O}) > d  & \Longleftrightarrow & (K) \lneqq (H) \, , \\
   c <  f(\mathcal{O}) < d & \Longleftrightarrow & (K) = (H) \, .
\end{eqnarray*}
In these equivalences, $(K)$ is the isotropy type of $\mathcal{O}$.
\end{lem}

\begin{proof} Assume that  $\mathcal{U}$ is a critical orbit 
of isotropy type $(L)\leq (H)$. Passage through $\mathcal{U}$ corresponds to 
attaching a $G$-handle of the form 
\[
   G \times_{L} (D(V) \times D(W)) \, . 
\]
The isotropy types occuring in this $G$-handle are smaller 
than or equal to $(L)$ and hence smaller than or equal to $(H)$. Now let $\mathcal{O}$ be a critical orbit 
and let $(K)$ be its isotropy type. 
Because we are dealing with special $G$-Morse functions, the attaching sphere
of the handle associated to $\mathcal{O}$  is of the form 
\[
   G \times_{K} (S^d \times 0)
\]
with $K$ acting trivially on $S^d$.  If $(K) \nleqq (H)$,  
then this attaching sphere must be disjoint 
from the first $G$-handle associated to $\mathcal{U}$. 
The order in which the two $G$-handles are attached 
can therefore be interchanged by 
adapting the given special $G$-Morse function appropriately. 
Hence we can assume that there is a noncritical value $c$ of $f$ such 
that for each critical orbit $\mathcal{O}$ of $f$, we have
\[
    f(\mathcal{O})  > c  \Longleftrightarrow  (K) \leq (H)  
\]
where $(K)$ is the isotropy type of $\mathcal{O}$. 
As before we can now argue that of the remaining $G$-handles, those 
of isotropy exactly $(H)$ are attached before those 
of isotropy type strictly smaller than $(H)$. This 
proves the existence of the second noncritical value $d$ with 
the stated property. 
\end{proof}

The next lemma gives an important  relation between the codimensions of $G$-handles and 
the codimensions of the associated $G$-strata. Again it is crucial 
to work with decompositions into special $G$-handles.

\begin{lem} \label{codestimate} Let $Z$ be a $G$-manifold 
with boundary and let 
\[  
    Z' =  Z  \cup_{\phi} \big( G \times_H (D^d  \times D(W)) \big)
\]
be obtained from $Z$ by attaching a special $G$-handle. 
Let $F \subset Z_{(H)}$ be the component containing 
$\phi(G \times_H (D^d \times 0))$. 
Then the codimension of $F$ in $Z$ and the dimension of $W$ 
are related by the inequality 
\[   
    {\rm codim}\,   F \leq \dim W \, . 
\]
\end{lem}

\begin{proof} Because the action of $H$ on 
$G \times  (D^{d} \times D(W))$ is free, we have   
\[
  \dim Z = \dim G + \dim  W + d  - \dim H \, . 
\]
But 
\[
   \dim G + d - \dim H \leq \dim F  
\]
becauce $H$ acts trivially on $D^{d}$ and hence $G \times_{H} 
( D^{d} \times 0)$ 
is completely contained in $F$. Assuming $\dim F \leq \dim Z -\dim W -1$, we 
therefore obtain
\[
    \dim Z \leq \dim F + \dim W \leq \dim Z -1  \, , 
\]
a contradiction. 
\end{proof} 

Finally we need an invariance statement for 
the $G$-homotopy type of certain singular strata 
under $G$-handle attachments.

\begin{lem} \label{attachnonmax} Let $Z$ be a $G$-manifold 
with boundary and let 
\[  
    Z' =  Z  \cup_{\phi} \big( G \times_H (D(V)  \times D(W)) \big)
\]
be obtained from $Z$ by attaching a $G$-handle such that 
$H$ acts trivially on $W$ (this condition is in some sense 
dual to that of being a special $G$-handle). 
If $(K) \neq (H)$, then there exists a $G$-homotopy 
equivalence
\[  
     Z_{(K)} \simeq Z'_{(K)} \, . 
\]
\end{lem}

\begin{proof} Because $H$  acts trivially on $W$ and 
$(K) \neq (H)$, we have
\[  
  \big( G \times_H (0 \times D(W)) \big)_{(K)}  = \emptyset
\]
and this (together with the fact that $V$ and $W$ are linear $H$-spaces) implies that 
the inclusion 
\[
      \big( G \times_H (S(V) \times D(W)) \big)_{(K)} \hookrightarrow 
      \big( G \times_H (D(V) \times D(W)) \big)_{(K)} 
\]
is a $G$-deformation retract.  
\end{proof}

Now we can formulate a bordism principle - Theorem A from the introduction -  for constructing 
$G$-invariant metrics of positive scalar curvature. In view of the 
Lawson-Yau theorem stated in the introduction, 
its main purpose is 
for Lie groups whose identity components are abelian. 
Recall that a singular stratum in a connected $G$-manifold $M$ is 
a connected component of some $M_{(H)}$ where $(H) \neq (H_{min})$. Singular 
strata are $G$-invariant submanifolds, but need  not be 
compact even if $M$ is compact (cf. Example \ref{typical}).  

\begin{thm} \label{possym} Let $Z$ be a compact connected oriented $G$-bordism
(with $G$ acting by orientation preserving maps) 
between the closed $G$-manifolds $X$ and $Y$. Assume the following: 
\begin{itemize}
      \item[i.)] ${\rm coh} (Z,G) \geq 6$ , 
      \item[ii.)] the inclusion $Y_{max} \hookrightarrow Z_{max}$ of 
               maximal orbits is a (nonequivariant) $2$-equivalence,   
      \item[iii.)] each singular stratum of codimension $2$ in $Z$ meets $Y$.
\end{itemize}
Then, if $X$ admits a $G$-invariant metric of positive scalar
curvature, the same is true for $Y$. 
\end{thm}

\begin{proof} Let  
\[
   f : Z \to [0,1]
\]
be a special $G$-Morse function. We will replace $f$ by a special $G$-Morse function 
without critical orbits of coindex $0$, $1$ or $2$. 
By Lemma \ref{codestimate}, all critical orbits of $f$ which 
are of this form are contained in singular strata of 
codimension less than $3$ in $Z$. Because $Z$ is oriented and $G$ acts in an 
orientation preserving fashion, there are no codimension-$1$ singular 
strata in $Z$. 

In a first step, we will take care of the singular strata of 
codimension $2$ that contain critical orbits of coindex less than $3$. 
Because $f$ is special, the coindex of these critical orbits 
is exactly $2$. 

Let  $F \subset Z$ be a singular stratum of codimension 
$2$ (which need not be compact) of isotropy type $(H)$ and containing  
coindex-$2$ critical orbits of $f$.  We will remove 
the critical orbits of coindex $2$ from $F$ without changing 
$f$ around the other singular strata in $Z$. 

By applying Lemma \ref{inducedhandledec} to the subgroup $H$, 
we can assume without loss of generality that there are  
noncritical values $0 < c < d < 1$ of $f$ such 
that for each critical orbit $\mathcal{O}$ of $f$ the equivalences 
\begin{eqnarray*}
     f(\mathcal{O}) > d & \Longleftrightarrow  & (K) \lneqq (H)      \, , \\ 
     c < f(\mathcal{O}) < d & \Longleftrightarrow & (K) = (H) 
\end{eqnarray*} 
hold where $(K)$ is the isotropy type of $\mathcal{O}$. Because 
$F$ is connected and $F \cap Y \neq \emptyset$, each component of 
$F \cap f^{-1}[0,d]$
has nonempty intersection with $f^{-1}(d)$. Furthermore, the map 
\[
    (f^{-1}[c,d])_{(H)} \hookrightarrow  (f^{-1}[0,d])_{(H)}
\]
is a $G$-homotopy equivalence by Lemma \ref{attachnonmax} (recall that 
$f$ is special, so $-f$ induces a decomposition into 
$G$-handles of the form described in Lemma \ref{attachnonmax}) and hence    
induces a bijection of connected components. This implies that each 
component $F_1 , \ldots, F_k$ of $F \cap f^{-1}[c,d]$ has nonempty intersection with $f^{-1}(d)$. 
We will now concentrate on the partial bordism 
\[
   P  : =   f^{-1}[c,d] \subset Z \, . 
\]
By construction, the critical orbits of the special $G$-Morse function 
\[
   f|_P : P \to [c,d]
\]
are exactly those critical orbits of $f$ which are of isotropy type $(H)$. 
Hence the restriction $f|_{P_{(H)}}$ induces a finite $G$-handle 
decomposition of $P_{(H)}$ and all critical orbits of $f|_F$ 
are contained in $P_{(H)}$. 
It is therefore enough to remove all coindex-$2$ critical
orbits of  $f|_P$ which are contained in $F \cap P$, but without changing $f|_P$ near the boundary
\[
   \partial P := f^{-1}(c) \cup f^{-1}(d) \subset P  
\]
and around singular strata in $P$ that are different from $F_1, \ldots, F_k$. 
We do this seperately for each component  $C \subset F \cap P$.

Using the fact that $C$ has nonempty intersection with $f^{-1}(d)$ (see above), 
we can use a (relative form of a) nonequivariant handle cancelation 
on the induced handle decomposition of $C/G$ (see e.g.~\cite{Milnorh}, Theorem 8.1)
in order to obtain a $G$-Morse function 
\[
   h  : C  \to [c,d]
\]
without coindex $0$ critical orbits, which coincides with $f|_C$ near 
$C \cap \partial P$ and outside a compact subset $K \subset C$. Note
that $C$ contains just one orbit type so that $h$ is 
automatically a special 
$G$-Morse function.
By (a relative form of) Lemma \ref{specialdense} we 
find a special $G$-Morse function 
\[
   P \to [c,d]
\]
whose restriction to $C$
coincides with $h$ and which is equal to $f|_P$ near $\partial P$ 
and near the singular strata in $P$ that are different from $C$. Because 
the codimension of $C$ in $P$ is $2$, the 
critical orbits in $C$ of this special $G$-Morse function are of coindex 
at least $1+2=3$.

This new special $G$-Morse function might have more critical orbits 
than $f$, but these do not lie on 
singular strata and are hence of minimal isotropy type. But 
critical orbits of minimal isotropy type (i.e.~lying in $Z_{max}$) 
will be taken care of later. 

We carry out the same process 
for all other components of $F \cap P$ and obtain a special 
$G$-Morse function $P \to [c,d]$ which coincides with 
$f|_P$ near $\partial P$ and near those singular strata
in $P$ that are different from a component of $F \cap P$. 
We combine this  new special 
$G$-Morse function on $P$ with the old $G$-Morse function $f|_{Z \setminus P}$ 
in order to produce a special $G$-Morse function 
\[
     Z \to [0,1]
\]
with no critical orbits of coindex $2$ in $F$,  which
coincides with $f$ near the 
singular strata different from $F$. 

After applying this procedure several times, we
get a special $G$-Morse function $f : Z \to [0,1]$ such that no  
singular stratum 
of codimension $2$ contains a critical orbit of coindex less than $3$. 

We now remove the critical orbits of coindex $0$, $1$ or $2$ in $Z_{max}$.
Arguing similarly as before we can assume that there is 
a noncritical value $c$ of $f$ so that  
\[
   P := f^{-1}[c,1]
\]
contains exactly those critical orbits of $f$ 
that are maximal. In particular, we get an induced 
handle decomposition of $P_{max}$. Because the inclusion 
\[
    Y_{max} \hookrightarrow P_{max}
\]
is a $2$-equivalence (this uses assumption ii.) in Theorem \ref{possym} as
well as Lemma \ref{attachnonmax}), the same is true for the 
inclusion of orbit spaces 
\[
    Y_{max}/G \hookrightarrow P_{max}/G  
\]
by comparing the long exact homotopy sequences 
for the respective $G/(H_{min})$-fibrations. 
We can now use nonequivariant handle cancelation on $P_{max}/G$ 
to find a handle decomposition of this space without 
codimension $0$, $1$ or $2$ handles. For the non-simply connected 
case, this is explained carefully in \cite{Lueck}. Note 
that by assumption, $P_{max}/G$ is oriented and of dimension at 
least $6$ so that the requirements for performing handle cancelation are 
fulfilled.

We end up with a special Morse function
\[   
   f : Z \to [0,1]
\]
without critical orbits of coindex less than $3$. 
The equivariant surgery principle (Theorem \ref{equivgl})
finishes the proof of Theorem \ref{possym}. 
\end{proof}

The connectivity assumption ii.) of Theorem \ref{possym}
reduces to a corresponding 
assumption in the nonequivariant bordism principle if 
$G$ is trivial.   
It is clear that the map $Y_{max} \hookrightarrow Z_{max}$ is a $2$-equivalence
if the inclusion $Y \hookrightarrow Z$ is a $2$-equivalence and $Z$ does 
not contain singular strata of codimension $2$ or $3$. 
However, we have not been able to replace condition ii.) by a connectivity 
assumption
on $Y \hookrightarrow Z$ if singular strata of codimension $2$ or $3$ occur in $Z$. 
Assumption iii.) can be viewed as the condition that the inclusion $Y \hookrightarrow 
Z$ restricted to singular strata of codimension $2$ be a $0$-equivalence (i.e. 
it induces a surjective map on $\pi_0$). 
We will make a few comments on the cohomogeneity restriction i.) of Theorem \ref{possym} 
at the end of this section. 

The following corollary of Theorem \ref{possym} is immediate. 

\begin{cor}[cf.~\cite{RosWei}, Theorem 2.3] \label{zetpe} Let $\Z/p$ 
act smoothly on a closed simply connected spin manifold $M^n$ 
where $n \geq 5$ and $p$ is an odd prime, preserving a spin structure. Assume 
furthermore that $M$ is equivariantly cobordant to another (not 
necessarily connected) spin $\Z/p$-manifold $M'$ by a bordism 
$W$ whose fixed set $W^{\Z/p}$ does not contain 
components of codimension $2$. 
If $M'$ has an invariant metric of positive scalar 
curvature, then so does $M$. 
\end{cor}

\begin{proof} Because  $p$ is odd, $W$ only contains fixed components
of even codimension and so the singular strata in $W$ have 
codimension  at least $4$. Because $\Z/p$ is a $0$-dimensional Lie group, 
the cohomogeneity assumption in Theorem \ref{possym} is satisfied. 
We can kill $\pi_1(W)$ and $\pi_2(W)$ by performing surgeries 
on the free part of the interior of $W$ and therefore assume 
that $M \hookrightarrow W$ is a $2$-equivalence. This
is equivalent to $M_{max} \hookrightarrow W_{max}$
being a $2$-equivalence because the singular strata in $W$ are 
of codimension at least $4$. 
\end{proof}

The original formulation in \cite{RosWei} requires that 
only the codimensions of singular strata contained in  $M$ be 
at least $4$.  The proof presented in 
{\it loc. cit.} argues  that the fixed set $W^{\Z/p}$ 
can be built from $(M')^{\Z/p} \times [0,1]$ by successive 
handle attachments.  These  handles are then  
thickened inside $W$ and thus replaced by handles 
of codimension at least $4$ in $W$ if the codimension 
of $W^{\Z/p}$ in $W$ (!) is assumed to be at least $4$.  
The remaining (free) part 
of $W$ is then constructed by successive handle attachments 
to the union of $M'$ and the thickening of $W^{\Z/p}$ obtained 
before.  The codimensions of these handles can be controlled 
by the topological assumptions on $M$.  However, we do not 
understand how this proof works if there are codimension-$2$ components 
in $W^{\Z/p}$ that are disjoint from $M$ (see the following 
Proposition \ref{illustrate}).  If we assume that each codimension-$2$ component in $W$ touches 
$M$ and 
\[
   M \setminus M^{\Z/p} \hookrightarrow W \setminus W^{\Z/p}
\]
is a $2$-equivalence, then our Theorem \ref{possym} shows that an invariant 
metric of positive scalar curvature on $M'$ can still be pushed 
through the bordism in order to produce one on $M$. So, 
in fact, our codimension-$2$ restriction is rather 
the opposite of the one  proposed in \cite{RosWei}, Theorem 2.3. 

The following proposition shows  that in the case of 
codimension-$2$ fixed components 
in $Z$ that do not touch $Y$, any special $\Z/p$-handle decomposition of $Z$  
contains handles of codimension $2$. 
This points towards a clear limitation of using (conventional) equivariant handle 
decompositions for constructing equivariant positive 
scalar curvature metrics and shows that a new idea 
is needed in order to deal with assumption iii.) in Theorem \ref{possym}.

\begin{prop} \label{illustrate} 
Let $Z$ be a compact $\Z/p$-bordism ($p$ an odd prime) between the closed $\Z/p$-manifolds
$X$ and $Y$. Assume that $Z$ contains a fixed component of codimension 
$2$ that is disjoint from $Y$. Then any 
$\Z/p$-handle decomposition of $Z$ (starting from $X \times [0,1]$)
contains handles of codimension $0$ or codimension $2$.
\end{prop}

\begin{proof} Each such handle decomposition is associated  to a $\Z/p$-Morse
function 
\[
    f : Z \to \R \, .  
\]
By assumption, the restriction of $f$ to the fixed point 
set $Z^{\Z/p}$ must have a local maximum $x \in F$ 
where $F \subset Z^{\Z/p}$ is a fixed component of codimension $2$
in $Z$ which is disjoint from  $Y$. Then $x$ 
is a critical orbit of $f$ of coindex $0$ or $2$ (it 
is of coindex exactly $2$ if $f$ happens to be special). 
\end{proof}

One might speculate that  Theorem \ref{possym} could be used to 
give an alternative proof of the Lawson-Yau theorem stated in 
the introduction. However, 
we have not been able to carry this out.  We remark again that in the proof 
of the Lawson-Yau theorem, an explicit metric 
with the required properties 
is constructed,  whereas the metric on $Y$ prescribed by  Theorem \ref{possym} 
depends on the given bordism $Z$ and is therefore difficult to make explicit.

We conclude this section with some remarks  on the 
cohomogeneity assumption i.) in Theorem 
\ref{possym}. The necessity of this assumption is obvious if $G$ is 
finite because a similar dimension restriction appears
in the nonequivariant case. We will therefore  
concentrate on  $S^1$-actions.  

Let $M$ be a closed  symplectic 
non-spin simply connected $4$-manifold 
with $b_2^+(M) \geq 2$. Then $M$ does 
not admit a positive scalar curvature metric by Seiberg-Witten theory
(see \cite{Nico}, Corollary 2.3.8 in combination with Theorem 3.3.29).
By \cite{Bebe}, Theorem C (cf.~Lemma  \ref{basic} below),
the free $S^1$-manifold $Y : = M \times S^1$ 
(with $S^1$ acting trivially on $M$) does not admit an 
$S^1$-invariant metric of positive scalar curvature, because
the quotient manifold 
\[
      (M \times S^1)/S^1 = M
\]
does not admit 
such a metric.  Moreover,  
the oriented bordism group $\Omega^{\SO}_4$ is 
generated by the bordism class of $\C P^2$ and hence 
the manifold $M$ is 
oriented bordant to a 
manifold $X$ which carries a metric of positive
scalar curvature. The inclusion $M \hookrightarrow W$ into such 
a bordism can be assumed to be a $2$-equivalence by performing 
appropriate surgeries on $W$, since $M$ is simply connected,
is not spin and $\dim W = 5$. The 
manifold $W \times S^1$ is then an $S^1$-bordism between 
$X \times S^1$ and $M \times S^1$, and moreover the $S^1$-manifold $X \times S^1$
admits an invariant metric of positive scalar curvature and the 
inclusion 
\[ 
   M \times S^1 \hookrightarrow W \times S^1
\]
is a $2$-equivalence. However, the cohomogeneity of $W \times S^1$ 
is equal to $5$ (and in particular smaller than $6$).

\section{Resolution of singularities} \label{genial}

We now specialize our discussion to fixed point free $S^1$-manifolds.
It turns out that in this case, the assumption on codimension-$2$ singular 
strata 
in Theorem \ref{possym} can be dropped if $Y_{max}$ is simply connected
and not spin. This improvement  is based on a new surgery technique 
which enables us to remove codimension-$2$ singular strata
from fixed point free $S^1$-manifolds while preserving  
a given invariant positive scalar curvature
metric. In this section, we will work out the details 
of this procedure before we discuss the improvement of Theorem \ref{possym} in 
the next Section \ref{anwend}. 

We will first describe the surgery step and discuss some
technical conditions needed in the geometric analysis 
of this situation before we finally state and prove the preservation 
of invariant positive scalar curvature metrics under 
the surgery in Theorem \ref{heikel}.   

Let $M$ be a closed fixed point free $S^1$-manifold 
of dimension $n \geq 3$ and let 
\[   
   \phi : S^1 \times_H (S^{n-3} \times D(W)) \hookrightarrow M 
\]
be an $S^1$-equivariant embedding where $H \subset S^1$ is a finite 
subgroup and $W$ is a one dimensional unitary effective $H$-representation. 
The group $H$ acts trivially on $S^{n-3}$. 

Because $S(W)/H$ can be identified with $S^1$, the $S^1$-principal 
bundle 
\[
    S^1 \hookrightarrow S^1 \times_H (S^{n-3} \times S(W)) \to S^{n-3} \times S(W)/H
\]
is trivial. After a choice of trivialization
\[
  \chi :  S^1 \times_H (S^{n-3} \times S(W)) \cong S^1 \times S^{n-3} \times S(W)/H
\]
and considering $S(W)/H$ as the boundary of $D^2$, we can glue 
the free $S^1$-manifold $S^1 \times S^{n-3} \times D^2$ back to 
$M \setminus \im (\phi)$ to get a new $S^1$-manifold $M'$. The manifold $M'$ 
is constructed from $M$ by removing the singular stratum $F := \phi(S^1 \times_H (S^{n-3} \times 0))$ via a kind of codimension-$2$
surgery.This surgery step is 
different from the equivariant surgery described at the 
beginning of Section \ref{equivgl}. 
We will refer to $M'$ as being obtained 
from $M$ by a {\em resolution} of the codimension-$2$ singular stratum $F$. We 
remark that in general this surgery step depends on the choice of the trivialization  $\chi$. 

Before proceeding, we need to impose some 
restrictions on the $S^1$-actions under consideration. 

\begin{defn}[cf. \cite{Ossa}, p. 46] \label{conditionce} Let $M$ be a compact $S^1$-manifold. We say 
that $M$ satisfies {\em condition $\C$} if for each closed subgroup 
$H \subset S^1$, 
the $S^1$-equivariant normal bundle of the closed  submanifold $M^H \subset M$ 
(which may contain different isotropy types) 
is equipped with the structure of a complex $S^1$-bundle such that the 
following compatibility condition holds: If  $K, H \subset S^1$ are 
two closed subgroups and $K \subset H$,  then the restriction of the normal bundle of $M^K \subset M$ 
to $M^H$ is a direct summand of the normal bundle of $M^H \subset M$ 
as a complex $S^1$-bundle. 
\end{defn}

Note that the singular strata of an $S^1$-manifold  
satisfying condition $\C$ always have even codimension. The following 
lemma shows that $S^1$-actions satisfying condition $\C$ occur 
naturally in many situations. 

\begin{lem} \label{oddorder} Let $M$ be a compact $S^1$-manifold for which all
finite isotropy groups have odd order. Then $M$ satisfies condition $\C$. 
\end{lem}

\begin{proof} Let $H \subset S^1$ be a closed subgroup and let $\nu \to M^H$ be 
the normal bundle of $M^H$ in $M$. For each $x \in M^H$, the fibre 
$\nu_x$ is a real $H$-representation with $\nu_x^H = 0$. 
This $H$-representation has a decomposition 
\[
  \nu_x = (E_1 \otimes V_1) \oplus \ldots \oplus (E_k \otimes V_k)
\]
where $E_i$, $1 \leq i \leq k$ are real vector spaces with trivial $H$-action 
and $V_i$ are pairwise different nontrivial
irreducible real $H$-representations. Because $H=S^1$ or $H$ is of odd order,
each $V_i$ has real dimension $2$  with 
$H$ acting as a rotation action. Hence,  each $V_i$ carries 
the structure of a one dimensional complex $H$-representation with 
$g \in H \subset S^1$ acting by 
\[   
    (g,v) \mapsto g^{\xi_i} \cdot v \, .   
\]
Here $\xi_i \in \Z$ is different from $0$ and not a multiple of $|H|$ if $H$ 
is a finite group. 
The complex structure on $V_i$ is uniquely determined if we require that 
$\xi_i > 0$ if $H = S^1$, resp. $\xi_i = 1, \ldots, \frac{|H|-1}{2} \mod |H|$ if 
$H$ is finite. We conclude that $\nu_x$ carries a canonical induced structure 
of a complex $H$-representation.  
\end{proof}

\begin{lem} \label{condcedec} Let $Z$ be a compact $S^1$-bordism 
satisfying condition $\C$ and assume that $Z$ 
is decomposed into special $S^1$-handles of the form  
\[
    S^1 \times_H (D^d \times D(W)) \, . 
\]
Then each of these $S^1$-handles is equipped with the 
following additional structure: If $K \subset S^1$ is a closed subgroup 
with $K \subset H$ and we decompose the orthogonal $H$-representation $W$ as 
\[
    W = W^K \oplus (W^K)^{\perp} \, , 
\]
then the $H$-representation $(W^K)^{\perp}$, which is $H$-invariant 
because $S^1$ is abelian, carries 
the structure of a unitary $H$-representation that is  
compatible with the given orthogonal $H$-structure. 
\end{lem}

\begin{proof}  This holds because each such $S^1$-handle is equivariantly 
embedded in $Z$ and the normal bundle of the $K$-fixed set
\[
 S^1 \times_H (D^d \times D(W)^K) = (S^1 \times_H (D^d \times D(W))^K \subset S^1 \times_H (D^d \times D(W))
\]
has fibre $(W^K)^{\perp}$. 
\end{proof} 

The next definition is of a rather
technical nature. These  properties of $S^1$-actions play a crucial role in our proofs, but do not appear in the final theorems.

\begin{defn} \label{technisch} Let $M$ be a manifold 
equipped with an $S^1$-action $\tau$ and with an $S^1$-invariant 
Riemannian metric $g$. 
\begin{itemize}
   \item[i.)] The metric is called {\em scaled} if the $S^1$-action 
             is fixed point free and the vector field on 
             $M$ generated by the action has constant length, called 
             the {\em scale} of the action.  
   \item[ii.)] We call $g$ {\em normally symmetric in codimension $2$} if 
              the following holds: Let $H \subset S^1$ be a closed subgroup  
              and let 
              $F \subset M^H$ be a component of the $H$-fixed 
              subset in $M$ which is of codimension $2$ in 
              $M$ (recall that $F$ is a 
              closed $S^1$-invariant 
              submanifold of $M$). Then there exists 
              an $S^1$-invariant tubular neighbourhood $N_F \subset M$ 
              of $F$ together with a second isometric 
              $S^1$-action $\sigma_F$ on $N_F$ that commutes with $\tau$ and 
              has fixed point set $F$.   
\end{itemize}
\end{defn}

If an $S^1$-invariant metric on $M$ has scale $s$, then the length of an orbit with isotropy $H$ 
is exactly $(2\pi s)/|H|$. We remind the reader of the following special case of the O'Neill 
formula (see \cite{Besse}, 9.37): Let 
\[
 S^1 \hookrightarrow   E  \stackrel{\pi}{\longrightarrow} B
\]
be a Riemannian submersion with totally geodesic fibres $S^1$ (necessarily of 
constant length, if $B$ is connected). Then 
the scalar curvatures of $E$ and $B$ are related by the formula 
\[
    s_E(x) = s_B(\pi(x)) - \|A(x)\|^2  
\]
where $x \in E$ and $A$ is the tensor field on $E$ defined in \cite{Besse}, 9.20. 
This makes clear why we prefer to work with invariant metrics 
of positive scalar curvature which are scaled: If a free $S^1$-manifold $E$
of dimension at least $3$ is equipped with such a metric, then the 
orbits are totally 
geodesic submanifolds of $E$ (see \cite{Besse}, Theorem 9.59)  and the induced metric 
on the orbit space $E/S^1$ has positive scalar curvature. This observation 
will be used at various places in the following discussion. 
In general, if the vertical part (tangent to the fibres) of the 
metric on $E$ is multiplied by a constant factor $\epsilon^2$ (we call the resulting 
Riemannian manifold $E_{\epsilon}$), 
then the norm $\|A\|^2$ is multiplied by $\epsilon^2$. 
In particular, shrinking the fibres $S^1$ increases the 
scalar curvature on $E$ and if $s_B > 0$ and $B$ is 
compact, then there is some constant $\epsilon_0 > 0$ so that 
$s_{E_{\epsilon}} > 0$, if $0 < \epsilon < \epsilon_0$. Also,  
if an $S^1$-invariant metric of positive scalar curvature 
is scaled, then the scale can be decreased arbitrarily while preserving 
the positive scalar curvature property. 

Note that the additional normal symmetries described in part ii.) 
of Definition \ref{technisch} need only be defined locally around the respective 
fixed point sets. These additional 
symmetries will 
considerably faciliate our later analytic arguments.  We consider normally 
symmetric $S^1$-invariant metrics as 
the ``naturally occuring'' ones. This will become clear in the proof 
of Lemma \ref{substantial2}. 
As a first illustration we provide the following example.

\begin{ex} \label{oversym} Let $V$ be a unitary $S^1$-representation.
Then the induced $S^1$-invariant metric on the 
representation sphere $S(V)$ (equipped with the restricted $S^1$-action)
is normally symmetric in codimension $2$. 
\end{ex} 

Neither of the restrictions formulated in Definition \ref{technisch} 
are serious. This is the content of the following two lemmas.  

\begin{lem} \label{substantial1} Let $M$ be a closed fixed point free $S^1$-manifold 
of dimension at least $3$ which is  equipped with an  
invariant metric of positive scalar curvture. Then $M$ 
also carries an invariant metric of positive scalar curvature which is scaled.
If the original metric is normally symmetric
in codimension $2$,  the same can be assumed for the new metric.  
\end{lem}

\begin{proof}  Let  
\[   
   X : M \to TM  
\]
be the vector field generated by the $S^1$-action. Because the action 
is fixed point free, $X$ has no zeros and  defines a $1$-dimensional 
subbundle $\mathcal{V} \subset TM$. Let the smooth 
function $f:M \to \R$ be defined by 
\[
    f(p) :=  \|X(p)\|_g \, . 
\]
We split $TM$ into $\mathcal{V}$ and its orthogonal complement $\mathcal{H}$. For
a moment we restrict attention to a tube 
\[
   S^1 \times_H D(V) \subset M 
\]
of the action. After pulling back the metric $g$ along the orbit map 
\[
   S^1 \times D(V) \to S^1 \times_H D(V) \, ,  
\]
we get a metric on $S^1 \times D(V)$ which is 
invariant under the free $S^1$-action on the first
factor. Let $h$ be the induced quotient metric on $D(V)$. We
set $n := \dim M$. The argument 
from \cite{Bebe}, Section 9,  shows that  the metric 
\[
      f^{\frac{2}{n -2}} \cdot h
\]
on $D(V)$ has positive scalar curvature.  

After this preparation, let $dt^2$ be the 
metric on $\mathcal{V}$ with respect to which $X$ has 
constant length $1$. An application of the O'Neill 
formula together with the previous local argument shows 
that there is a constant $\epsilon_0 > 0$ so that the metric 
\[
     (\epsilon^2 \cdot dt^2) \oplus (f^{\frac{2}{n - 2}} \cdot g|_{\mathcal{H}}) 
\]
on $M$ has positive scalar curvature if $0 < \epsilon < \epsilon_0$. 
By construction, this new metric has all of the required properties. 
If $g$ is normally symmetric, this follows because the additional 
$S^1$-actions $\sigma_F$ around the different 
codimension-$2$ subsets $F\subset M$
respect the decomposition $TM = \mathcal{V} \oplus \mathcal{H}$.  
\end{proof}

Concerning the construction of normally symmetric metrics of positive scalar curvature, 
we have the following variant of the surgery principle from Section \ref{equivGL}. 
We formulate this result directly in the form needed in Section \ref{anwend}. 

\begin{lem} \label{substantial2} Let $Z$ be a compact $S^1$-bordism 
satisfying condition $\C$ between  the 
closed $S^1$-manifolds $X$ and  $Y$. Assume that $X$ carries 
an invariant metric of positive scalar curvature 
which is normally symmetric in codimension $2$. If $Z$ admits 
a decomposition into special $S^1$-handles (starting from 
$X \times [0,1]$) of codimension at
least $3$, then $Y$ also carries an invariant metric of 
positive scalar curvature which is normally symmetric in 
codimension $2$.
\end{lem}

\begin{proof}  The decomposition of $Z$ into special 
$S^1$-handles of codimension at least $3$ implies that $Y$ is obtained 
from $X$ by performing equivariant surgeries of codimension at least 
$3$ along embedded submanifolds of the form 
\[
   S^1 \times_H (S^{d-1} \times D(W))  \, .  
\]
If $d=0$, then  this attaching locus 
is empty and the surgery step produces a new 
component $S :=  S^1 \times_H (D^0 \times S(W))$.  
We show at  first that these $S^1$-manifolds $S$ admit $S^1$-invariant 
metrics of positive scalar curvature which are normally symmetric in 
codimension $2$. 

Let $K \subset S^1$ be closed and let $F = S^K \subset S$ be a 
component of codimension $2$. Then $K \subset H$ (otherwise 
$S^K$ would be empty) and $W$ admits an orthogonal splitting 
\[
   W = W^K \oplus (W^K)^{\perp} \, .  
\]
Because $Z$ satisfies condition $\C$, Lemma \ref{condcedec} implies that $(W^K)^{\perp}$
has the induced structure of a unitary $H$-representation. 

We now consider 
the orthogonal $S^1$-action on 
\[   
   W = W^K \oplus (W^K)^{\perp}
\]
which is given by complex multiplication on the second summand. The
restriction of this  $S^1$-action 
to $S(W)$ commutes with the $H$-action. Furthermore, this $S^1$-action on $S(W)$ has fixed point set 
$S(W)^K$ and leaves the usual round metric on $S(W)$ invariant. In this way 
we obtain an  $S^1$-action $\sigma$ on $S$ by letting $S^1$ act on the 
$S(W)$-factor of $S^1 \times_H (D^0 \times S(W))$. This 
action commutes with the original $S^1$ action on $S$ (acting on the $S^1$-factor),
leaves the positive scalar curvature 
metric on $S$ constructed in the proof of Theorem \ref{equivgl} invariant (recall that this is 
a Riemannian submersion metric on the fibre bundle 
$S(W) \hookrightarrow S^1 \times_H (D^0 \times S(W)) \to S^1/H$)  and has 
fixed point set 
\[
    S^1 \times_H (D^0 \times S(W^K)) =  S^K \, . 
\]
This finishes the discussion of $S^1$-handles of dimension $0$ in $Z$. 

From now on, we concentrate on $S^1$-handles in $Z$  that are of dimension at least $1$. 

We need to prove the following fact: 
Let $M$ be a closed $S^1$-manifold equipped with an invariant 
metric of positive scalar curvature which is normally
symmetric in codimension $2$. If $M'$ is  obtained from 
$M$ by equivariant surgery of codimension at least $3$ along an 
$S^1$-equivariant embedding 
\[
   \phi :  S^1 \times_H (S^{d-1} \times D(W)) \hookrightarrow M  
\]
with $d-1 \geq 0$ (i.e. $S^{d-1} \neq \emptyset$), then $M'$ also carries an invariant metric of positive
scalar curvature which is normally symmetric in codimension $2$.  We argue as follows. 

At first notice that it suffices to treat the case when  there is a closed subgroup $K \subset S^1$ and a 
codimension-$2$ 
component $ F \subset M^K$ (for brevity we will call such submanifolds {\em fixed codimension-$2$ components}) such that 
\[
    \phi(S^1 \times_H (S^{d-1} \times 0)) \cap F \neq \emptyset \, . 
\]
Otherwise we 
could assume that $\im (\phi)$ would be  disjoint from any 
fixed codimension-$2$ component of $M$ and this would imply 
that the $S^1$-handle $S^1 \times_H (D^{d} \times D(W))$ would have no 
fixed codimension-$2$ components, either. Hence there would be nothing 
to prove. 

We may assume (up to $G$-isotopy) that 
$\phi$ maps the radial lines in the fibres of the $S^1$-equivariant fibre 
bundle  
\[
   D(W) \hookrightarrow S^1 \times_H (S^{d-1} \times D(W)) \to S^1 \times_H 
  (S^{d-1} \times 0)
\]
to unit speed geodesics in $M$ that are orthogonal to $\phi(S^1 \times_H (S^{d-1} \times 0))$. 
Now let 
\[   
   F_i \subset M^{K_i} \, , i = 1, \ldots, r\, , 
\]
be those fixed codimension-$2$ components (with certain -
necessarily distinct - subgroups 
$K_i \subset S^1$) 
which intersect $\phi(S^1 \times_H (S^{d-1} \times 0))$ nontrivially. This implies 
that  $K_i \subset H$ for all $i$ and 
\[
    \phi(S^1 \times_H (S^{d-1} \times 0)) \subset \bigcap_i F_i  \, . 
\]
Without loss of generality we may assume that 
\[
   \im (\phi) \subset N_{F_i}
\]
for all $i = 1, \ldots, r$, where $N_{F_i}$ are tubular neighbourhoods of $F_i$ which 
are equipped with additional actions $\sigma_i$ as described in Definition \ref{technisch}. 

The map $\phi$ and the actions $\sigma_1, \ldots, \sigma_r$ 
induce $S^1$-actions $\Sigma_i$ on $S^1 \times_H (S^{d-1} \times D(W))$ 
which are induced by orthogonal $S^1$-actions on $W$ (with respect to the 
standard Euclidean metric on $W$)
commuting with the $H$-action on $W$ and 
with fixed points sets $W^{K_i}$.
Let $g$ be the given metric on $M$ and 
$\phi^*(g)$ be the induced metric on $S^1 \times_H (S^{d-1} \times D(W))$. 
This metric enjoys the following $S^1$-symmetries: 
\begin{itemize}
    \item[i.)] It is invariant under the  $S^1$-action on the first factor of
              $S^1 \times_H (S^{d-1} \times D(W))$,
    \item[ii.)] for each $i=1, \ldots, r$ it is invariant under the action $\Sigma_i$. 
\end{itemize}

The constructions in the proof of Theorem \ref{equivgl} 
preserve all these $S^1$-symmetries of  $\phi^*(g)$: 
This is clear for the submersion metrics on the sphere bundles
\[
   S_{\epsilon}(W) \hookrightarrow S^1 \times_H (S^{d-1} \times S_{\epsilon}(W)) 
   \to S^1 \times_H (S^{d-1} \times 0)
\]
which are constructed in the proof of Theorem \ref{equivgl} and hence also 
for the metric on 
\[
    S^1 \times_H (S^{d-1} \times S(W)) \times [0,T]
\]
which interpolates (via the bending outward process) between the metric 
$\phi^*(g)$ and a submersion metric of the above form (with small $\epsilon$).

This finishes the proof of Lemma \ref{substantial2}
\end{proof}

The main  result of this section reads as follows. 

\begin{thm} \label{heikel} Let $M^n$ be a closed fixed point free $S^1$-manifold 
of dimension $n \geq 3$ and let 
\[
   \phi : S^1 \times_H (S^{n-3} \times D(W)) \hookrightarrow M 
\]
be an $S^1$-equivariant embedding where $W$ is a unitary effective $H$-representation 
of dimension $1$.  Let the manifold $M'$ be obtained from $M$ by 
resolving the singular stratum $\phi(S^1 \times_H (S^{n-3} \times 0)) 
\subset M$. 
If the manifold $M$ admits an invariant metric of positive scalar curvature
which is scaled and normally symmetric in codimension $2$, then also 
$M'$ admits such a metric.  
\end{thm} 

The remainder of this section is devoted to a proof of Theorem \ref{heikel}. 

Let $g$ be an invariant metric of positive scalar 
curvature on $M$ which is scaled and normally symmetric in codimension $2$. 
As before we assume that $\phi$ maps the 
radial lines in $D(W)$ to unit speed geodesics in $M$ 
orthogonal to $\phi(S^1 \times_H (S^{n-3} \times 0)) \subset M$ 
and that $\im (\phi)$ is contained in an $S^1$-invariant 
tubular neighbourhood of $\phi(S^1 \times_H (S^{n-3} \times 0))$ 
that is equipped with an additional $S^1$-action $\sigma$ as described
in Definition \ref{technisch}. We still denote the original $S^1$-action 
on $M$ by $\tau$. We thus obtain corresponding actions $\tau$ and $\sigma$ on 
\[
        N := S^1 \times_H (S^{n-3} \times D(W))
\] 
that are induced by rotation actions on the factors $S^1$ and $D(W)$ respectively. 
Because $W$ is a unitary $H$-representation, 
the induced $H$-action on $D(W)$ commutes with $\sigma$ (this also follows from 
the requirement that the actions $\tau$ and $\sigma$ commute). 

Recall that the induced metric $\phi^*(g)$ (that we will denote by $g$ from now on) 
on the total space of the $S^1$-equivariant fibre bundle 
\[
    D(W) \hookrightarrow N  \to S^1 \times_H (S^{n-3} \times 0) 
\]
need not be a Riemannian submersion metric 
and - contrary to the surgery principle explained in Section \ref{equivGL} - 
the consideration of such a metric does not seem to be of much use for our purposes 
because the representation sphere $S(W)$ is of dimension one
and hence does not carry a Riemannian metric of positive scalar curvature. 

From now on, we will write $S$ instead of $S^{n-3}$ for simplicity. On the orbit space
\[
    N/(S^1, \tau) = S  \times D(W)/H  
\]
we obtain an induced metric of positive scalar curvature away from the 
singular locus $S \times 0$ because $g$ is scaled. In a first 
step we will deform this metric near the singular locus 
so that it can be extended to a smooth metric  of 
positive scalar curvature on $S \times D^2$. Here we identify the cone factor $D(W)/H$ with 
$D^2$. Roughly speaking, this deformation is possible because the 
tip of the cone $D(W)/H$ can be viewed as a source of a large amount 
of positive scalar curvature which can be distributed over 
a neighbourhood of the singularity in $D(W)/H$.

One of the main technical problems is the extension of 
the induced metric on $S \times (D(W) \setminus 0)/H$ to a smooth metric on 
$S \times D^2$ in a well controlled way. At this point we make 
essential use of the additional rotation symmetry (induced by $\sigma$)
of this metric.  

The metric $g$ can be pulled back along the quotient  
map 
\[
   S^1  \times (S  \times D(W)) \to S^1 \times_H  (S \times D(W))  
\]
and yields a metric on $S^1 \times (S  \times D(W))$ which 
is invariant under the usual $S^1$-rotation actions  on $S^1$ and on $D(W)$. 

Ater dividing out the (free) $S^1$-action on the $S^1$-factor 
(this action is lifted from $\tau$), we get an induced metric $h$ on 
$S \times D(W)$ which is 
of positive scalar curvature (because the action $\tau$ is scaled) and invariant under the 
rotation action on the $D(W)$-factor.

The space $N/(S^1,\tau)$ 
can now be identified with the orbit space of $S\times D(W)$ under 
the $H$-action on $D(W)$. We consider the (in general not orthogonal)
canonical splitting of bundles over $S \times D(W)$ 
\[
     T(S \times D(W)) \cong TS \oplus T(D(W)) \, . 
\]
Using local coordinates $u_1, \ldots, u_m$ (where $m = n-3$) on $S$ and polar coordinates
$(r,\theta)$ on $D(W) \setminus 0$, the restriction of $h$ to $S \times (D(W) \setminus 0)$ 
can be written as
\[
   \sum_{1 \leq i \leq  j \leq m} \alpha_{ij}(u,r) du_i du_j + \sum_{1 \leq i \leq m} \beta_{i}(u,r) 
   du_i d\theta + dr^2 + \gamma(u,r)^2  d\theta^2 \, .   
\]
Summands of the form $du_i dr$ are not needed here because 
the radial lines in $(S \times D(W), h)$
are orthogonal to $S \times S_{r}^1$ for each $0 < r \leq 1$. 
Here and in what follows, $S_{r}^1$ will be  the $r$-sphere 
in $D(W) = D^2 \subset \C$. Furthermore, because $h$ is invariant under rotation of $D(W)$, the coefficient 
$\gamma$ does not depend on $\theta$. For later use we note 

\begin{lem} \label{umrechnen} Each of the summands 
\[
   \sum_{1 \leq i \leq m} \beta_{i}(u,r) du_i d\theta
\]
and 
\[
     dr^2 + \gamma(u,r)^2  d\theta^2
\]
extends to  a smooth  $(0,2)$-tensor on $S \times D(W)$. 
\end{lem}

\begin{proof} Under the conversion of polar into cartesian coordinates 
\begin{eqnarray*} 
    x & = & r \cos \theta \, , \\
    y & = & r \sin \theta
\end{eqnarray*} 
on $D(W)$, the first of the displayed terms is transformed into 
the sum of those components of the metric $2$-tensor 
of $(S \times (D(W) \setminus 0), h)$  
containing $du_i dx$ or $du_i dy$ and 
the second of the displayed terms is transformed into the 
sum of those components containing $dx^2$, $dx dy$ and $dy^2$.
Because $h|_{S \times (D(W) \setminus 0)}$ obviously extends to $S \times D(W)$,
each of these seperate summands extends to a smooth $(0,2)$-tensor 
on $S \times D(W)$. 
\end{proof}

Using the standard identification $D(W)/H = D^2$ which is induced by the map 
\[
    S^1 \to S^1 \, , ~ x \mapsto x^{|H|} \, ,  
\]
(here we use effectiveness of the $H$-action on $W$) we get a homeomorphism 
\[
    S \times D(W)/H \approx S \times D^2 \, ,  
\]
which is a diffeomorphism on $S \times (D(W) \setminus 0)/H$. 
Hence (and because $W$ is a unitary $H$-representation),
the metric $h$ induces a metric $q$ on
$S \times (D^2 \setminus 0)$ which is given by  
\[    
   q = \sum  \alpha_{ij}  du_i du_j + \frac{1}{|H|} \sum \beta_{i} du_i d\theta + dr^2 + 
    \frac{\gamma^2}{|H|^2} d\theta^2
\]
using again polar coordinates $(r,\theta)$ on $D^2 \setminus 0$
and the fact that $\gamma$ is independent of $\theta$. 
This metric cannot be extended to a  smooth metric on $S \times D^2$. 
One reason is that the  partial derivative
\[
    \frac{1}{|H|} \cdot \frac{ \partial \gamma(u,r) }{\partial r}|_{r=0}  =  \frac{1}{|H|}
\]
is different from $1$ (this fact corresponds to the conelike form of this metric). 
For a general comparison of 
metrics given in polar and cartesian coordinates, 
see \cite{Pet}, 1.3.4. However, this failure can 
be remedied by using a ``bending inwards''-process of the following form. 

We fix the slope 
\[ 
     c := \frac{\sqrt{|H|^2-1}}{|H|}
\]
and consider the affine function 
\[
   \kappa :  [0,1] \to \R \, , ~ r \mapsto -c r + c  \, . 
\]
Taking into account that 
\[
    \frac{1}{|H|^2} + c^2 = 1 \, , 
\]
the metric $q$ away from the singular locus $S \times 0$ is the induced metric 
on the hypersurface 
\[  
     \{ ((u,r,\theta),t) \in (S \times (D^2 \setminus 0) , \widetilde{q} ) \times (\R,dt^2) ~|~
      (r,t) = (r, \kappa(r)) \} \, .  
\]
Here we use the metric 
\[
   \widetilde{q}  : =  \sum  \alpha_{ij} du_i du_j + \frac{1}{|H|} \sum_i \beta_{i} du_i d\theta 
+ \frac{1}{|H|^2}(dr^2 + \gamma^2 d\theta^2) 
\]
on $S \times (D^2 \setminus 0)$ and the canonical metric $dt^2$ on $\R$. 

The first statement of the following lemma is the main reason 
why $g$ was  assumed to be normally symmetric.

\begin{lem} \label{eigQ} The metric $\widetilde{q}$ extends to a smooth metric on $S \times D^2$.
With respect to this metric, the subset $S \times 0 \subset S \times D^2$ 
is a totally geodesic submanifold and the radial lines 
in $S \times D^2$ are geodesics of constant speed $\frac{1}{|H|}$ 
\end{lem}

\begin{proof} The first claim follows from Lemma \ref{umrechnen}. Moreover, 
$S \times 0$ is a totally geoedesic submanifold because  
$\tilde{q}$ is invariant under the rotation action on the $D^2$-factor
induced by $\sigma$.  

For proving the last claim, we must show that
\[
    \nabla_{\partial_r} \partial_r = 0   
\]
where $\nabla$ is the connection induced by $\tilde{q}$. 
But this follows directly from the Koszul formula, the definition of $\tilde{q}$ and 
the fact that radial lines are (unit speed) geodesics for the metric $q$.
\end{proof}

In this picture, the bending inwards process amounts to replacing the preceding 
hypersurface by a smooth Riemannian submanifod 
\[
   \Sigma \subset (S \times D^2, \widetilde{q}) \times (\R, dt^2) 
\]
as follows. Let 
\[
   \lambda: [0,1] \to \R_{\geq 0} \, , r \mapsto \lambda(r) \, , 
\]
be a smooth function which has the following properties: There is a positive 
constant $C$ so that   
\begin{itemize}
    \item[i.)] the function $[-1,1] \mapsto \R$, $s \mapsto \lambda(|s|)$, is smooth,  
    \item[ii.)] the second derivative of $\lambda$ is nonpositive and 
             bounded in absolute value by $2C$,       
    \item[iii.)] there are positive constants $\epsilon, \mu$ so that  
             $\lambda(r) = \kappa(r)$, if  $\mu + \epsilon \leq r \leq 1$ 
             and $\lambda(r) = \lambda(0) - Cr^2$ for $0 \leq r \leq \mu$. 
    \end{itemize}
The numbers $\mu$, $\epsilon$ and $C$ are to be specified later. 
For brevity we will refer to such a function 
$\lambda$ as a {\rm profile} of {\em width} $\mu$, with {\em bending parameter} $C$ and of 
{\em adjusting length} $\epsilon$. The following 
fact is elementary. 

\begin{lem} \label{profil} Let $C,\mu > 0$ be chosen such that 
\[
       2  C \cdot \mu > c \, . 
\]
Then for any $\epsilon > 0$, there exists a profile $\lambda$ with bending 
parameter $C$, of width smaller than $\mu$ and of 
adjusting length smaller than $\epsilon$. 
\end{lem}

For $r \in [0,1]$ we denote by 
$\Theta(r) \in [0,\frac{\pi}{2})$ the angle between the graph of $\lambda$ and 
the $r$-axis at the point $(r, \lambda(r))$ and by 
\[
    k(r) = \frac{d \Theta}{d \xi}
\]
the curvature of ${\rm graph} (\lambda)$ at $(r,\lambda(r))$. 
Here $\xi$ parametrizes the graph of $\lambda$ as
a unit speed curve in $\R^2$ with initial point $(1,0)$. 
Because $\lambda'' \leq 0$, the function $k(r)$ is nonpositive
for any profile $\lambda$. 

\begin{figure}[htb]
\begin{center}
\includegraphics[scale=0.9]{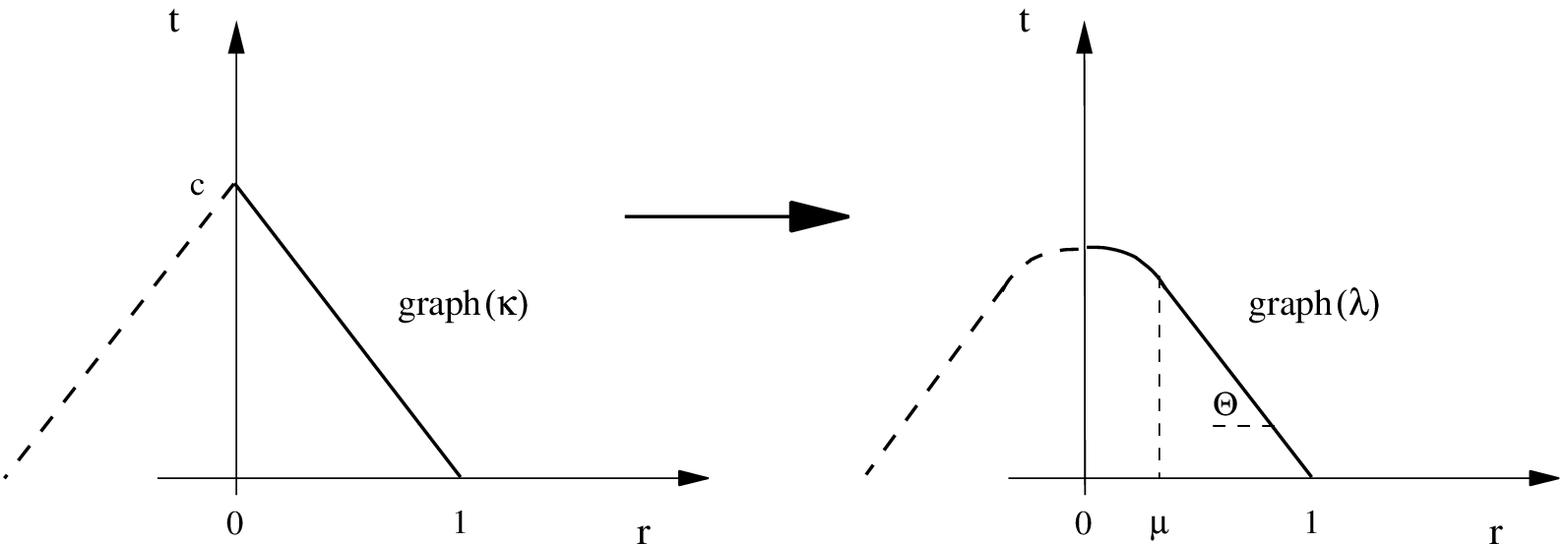}
\end{center}
\end{figure}

Depending on a given profile $\lambda$, we now consider the smooth
Riemannian submanifold 
\[
   \Sigma  := \{ ((u,r,\theta),t) \in (S \times D^2 , \widetilde{q}) \times (\R,dt^2) ~|~(r,t) \in
    {\rm graph}(\lambda) \} \, . 
\]
Let $U \subset S$ be an open subset and  
\[
   (e_1, \ldots, e_m) \, , ~ e_i \in \Gamma(TU) \, , 
\]
be an orthonormal frame with respect to
the metric $\widetilde{q}|_{S \times 0} (= h|_{S \times 0})$.
We get a family $(e_1, \ldots, e_m)$ of 
orthonormal fields on $U \times D^2$ that are tangential to the 
submanifolds $S \times S^1_r$, $r \in (0,1]$, by parallel transport along 
radial lines. The fields
\[
   e_{m+1}'(u,r,\theta) := \frac{|H|}{\gamma(u,r)} \cdot \partial_{\theta} 
   \, , ~~
   e_{m+2} : = |H| \cdot \partial_r
\]
complete this to a family of vector fields 
\[   
   (e_1, \ldots, e_m , e_{m+1}', e_{m+2})
\]
defined on $U \times (D^2\setminus 0)$ and of unit length. This need
not be an orthonormal frame. Indeed, $e_{m+1}'$ is  
orthogonal to $e_{m+2}$, but it need not be orthogonal to the 
fields $e_{1}, \ldots, e_m$, as the submanifolds $\{u \} \times D^2 \subset 
U \times D^2$ need not be totally geodesic. However, there is 
a vector field $e_{m+1}$ on $U \times (D^2 \setminus 0)$ that is 
tangential to the submanifolds $S \times S^1_r$, that  satisfies 
\[
     \lim_{r \to 0} \| e_{m+1}(u,r,\theta) - e_{m+1}'(u,r,\theta) \| = 0 
\]
for each fixed $(u,\theta)$ and that yields an orthonormal frame 
\[
   (e_1, \ldots, e_m, e_{m+1}, e_{m+2}) 
\]
on $U \times (D^2 \setminus 0)$: The vector 
$e_{m+1}(u,r,\theta)$ is given by parallel transport along 
the radial line 
\[
     l_{\theta} := \{ (u,r,\theta) ~|~ r \in [0,1] \} \subset U \times D^2
\]
of the vector in $T_{(u,0)}(U \times D^2)$ which is tangential 
to $\{u\} \times D^2$ (hence orthogonal to $S \times 0$) and includes 
an angle of $\pi/2$ with $l_{\theta}$.   
We denote by $s_{S \times D^2}$ 
the scalar curvature of $(S\times D^2, \widetilde{q})$ and by 
$s_{S \times S_r^1}$ the scalar curvature
of $(S \times S^1_r, \widetilde{q}|_{S \times S_r^1})$. Furthermore, 
for $r \in (0,1]$, let   
\[
   K_{S \times S^1_r}(u,\theta)  = \sum_{1 \leq i \leq m+1} \langle \nabla_{e_i} e_i, 
e_{m+2} \rangle 
\]
denote the mean curvature of the submanifold 
\[
    S \times S^1_r \subset (S \times D^2, \widetilde{q})   
\]
at $(u,\theta)$. With these specifications, the scalar curvature  
$s_{\Sigma}$ of $\Sigma$ at $(u,r,\theta)$ 
is given by the following formula. 

\begin{lem}[\cite{Wehrheim}, formula (4.1)] \label{jan} At the point 
$(u,r,\theta) \in U \times D^2$, we have
\[
 s_{\Sigma} = \cos^2 \Theta(r) \cdot s_{S \times D^2} + \sin^2 \Theta(r) \cdot
   s_{S \times S^1_r} + 2k(r) \sin \Theta(r)
 \cdot K_{S \times S^1_r}  \, . 
\]
\end{lem} 

We need to describe the behavior of the relevant geometric quantities near 
$r = 0$. 

\begin{lem} \label{asymptotic} \begin{itemize}
\item[i.)] For fixed $(u,\theta)$ and varying $r$, we have $K_{S \times S_r^1} 
          = - \frac{|H|}{r} + O(1)$. Furthermore, the remainder term $O(1)$ 
          depends continuously on $(u, \theta)$ and hence defines a continuous 
          function $S\times D^2 \to \R$ which vanishes on $S \times 0$.                                    
 \item[ii.)] The scalar curvature $s_{S \times S_r^1}$ 
              is  uniformly bounded with respect to  $r >  0$ and 
             all $(u,\theta)$.
\end{itemize}

\end{lem}

\begin{proof} For fixed $(u,\theta)$ and varying $r \in (0,1]$, we 
consider  the previously described orthonormal base $(e_1, \ldots, e_m, e_{m+1}, e_{m+2})$ of 
$T_{(u,r,\theta)} (S \times D^2)$. Recall  
that $e_1, \ldots, e_m$ are actually smooth vector fields on the whole 
of $U \times D^2$. For $1 \leq i,j \leq m+1$, let 
\[  
   \alpha_{ij}(r) : = \langle \nabla_{e_i} e_j , e_{m+2}  \rangle  
\]
be the components of the second fundamental form of $S \times S_r^1 \subset (S \times D^2, \widetilde{q})$ 
at the point $(u,r,\theta)$. 
In a first step we show that for  $1 \leq i,j \leq m$ we have
\begin{itemize}
    \item[i.)] $\alpha_{ij}(r) = O(r)$ ,   
    \item[ii.)] $\alpha_{m+1, j}(r) = O(1)$ ,   
    \item[iii.)] $\alpha_{m+1, m+1}(r) = - \frac{|H|}{r} + O(1)$ .   
\end{itemize}
These expansions are proved as follows. For each $1 \leq i,j \leq m$ the expression 
$\nabla_{e_i} e_j$ defines a smooth vector field on $U \times D^2$ 
which is tangential to $S \times 0$, because this is a totally geodesic 
submanifold of $S \times D^2$ (see Lemma \ref{eigQ}). This proves 
equation i.). 

For equation ii.) notice that the restriction 
\[
     e_{m+1}|_{U \times (0,1] \times \{ \theta \} } 
\]
can be extended to a smooth vector field on $U \times [-1,1] \times \{\theta\}$ 
by parallel transport along radial lines. Hence, the smooth function 
\[
   \alpha_{m+1, j}(u,-,\theta) : (0,1] \to \R 
\]
extends to $[-1,1]$. 
 
For proving the third equation, it is enough (using points i.) and ii.)) to show the asymptotic 
expansion 
\[   
  \langle \nabla_{e_{m+1}'} e_{m+1}', e_{m+2} \rangle = - \frac{|H|}{r} + O(1) 
\]
we recall (see Lemma \ref{eigQ}) that the 
radial lines 
\[
    [0,1] \to (S \times D^2, \widetilde{q}) \, , ~~ r \mapsto (u,r,\theta)
\]
are geodesics of constant speed $\frac{1}{|H|}$. 
We can therefore use the proof of  \cite{GL}, Lemma 1, applied to the submanifold 
$\{ u \} \times D^2 \subset S \times D^2$ with the restricted metric.   

It is clear, that the remainder terms in the expansions i.), ii.) and iii.) depend 
continuously on $(u,\theta)$.

Using these asymptotic expansions, the first claim of Lemma \ref{asymptotic} is immediate
and the second claim follows from the Gauss equation 
\[
 s_{S \times S^1_{r}} = 2 \cdot ( \sum_{1 \leq i < j \leq m+1} K(e_i,e_j)+ \alpha_{ii} \alpha_{jj} -
\alpha_{ij}^2) \, , 
\]
where $K(e_i,e_j)$ denotes the sectional curvature of the plane spanned by $e_i$ and $e_j$. 
Here it is crucial that  
\[
    \alpha_{ii}(r) \alpha_{m+1,m+1}(r) = O(1) 
\] 
for all $1 \leq i \leq m$ and that $\widetilde{q}$ is a smooth metric on the whole 
of $S \times D^2$ which is (of course) bounded in the  $C^{2}$-norm 
so that these sectional curvatures are uniformly bounded on $S \times D^2$. 
\end{proof}

We remark that in principle, this proof works for any 
smooth metric on $S \times D^2$
with respect to which $S \times 0$ is totally geodesic. The normal 
symmetry of $g$ was only used in order to construct such a metric from 
the conelike metric on $S \times (D(W)/H)$ induced by $g$.

We are now in a position to show that the profile $\lambda$ can be chosen such 
that the corresponding 
hypersurface $\Sigma$ has positive scalar curvature.
First we observe that for a profile $\lambda$ with bending parameter $C$, 
Lemma \ref{asymptotic} implies for all $(u,\theta)$ the important asymptotic 
expansion
\[
  2 k(r) \sin \Theta(r) \cdot K_{S \times S_r^1}(u,\theta) = (- 4C + O(r)) \cdot (2Cr+O(r^2)) \cdot \big( -\frac{|H|}{r} 
   +O(1)\big) =  8 |H| C^2 + O(r)  \, .  
\]
Furthermore, the function on the left hand side extends to a continuous function on 
$S \times D^2$.

\begin{lem} \label{uniform} Let $\lambda_1$, $\lambda_2$ be two profiles with bending 
parameters $C_1$ and $C_2$ and of widths $\mu_1$ and $\mu_2$, respectively. 
If $C_1 \leq C_2$,  then 
\[
   0 \geq  2k_1(r) \sin \Theta_1(r) \geq  2k_2(r) \sin \Theta_2(r) 
\]
for $0 \leq r \leq \min(\mu_1, \mu_2)$.  
\end{lem} 

\begin{proof} for $0 \leq r \leq \min(\mu_1, \mu_2)$, both profiles are 
in the standard form 
\[
    \lambda_1(r) = \lambda_1(0) - C_1 r^2 \, , ~ \lambda_2(r) = \lambda_2(0)  - C_2 r^2 \, . 
\]
\end{proof} 

This lemma together with the preceding asymptotic expansion make clear that there 
are $\overline{C}, \overline{\mu} > 0$ so that 
for any profile $\lambda$ with bending parameter $C \geq \overline{C}$ and width 
$\mu \leq  \overline{\mu}$, we have 
\[
   2k(r) \sin \Theta(r) \cdot K_{S \times S^1_r}(u,\theta) \geq |s_{S \times D^2}(u,r,\theta)|
 + |s_{S \times S^1_r}(u,\theta)| +1 
\]
for all $0 < r \leq \mu$ and all $(u, \theta)$. Note that the right hand side of this 
equation as well as the value of $K_{S \times S^1_r}(u,r,\theta)$                   
are independent of the particular profile $\lambda$. 
 Using Lemma \ref{profil} we can choose a profile with 
data $C, \mu$ satisfying these conditions. Furthermore,  we get 
\[
    K_{S \times S^1_r}(u,\theta) < 0 
\]
for $0 < r \leq \overline{\mu}$ and all $(u, \theta)$.

Again referring to Lemma \ref{profil},  
we can assume that the adjusting 
length of $\lambda$ is arbitrarily small. This is helpful because 
there is a constant $s> 0$ so 
that for any profile $\lambda$ we have 
\[
   \cos^2 \Theta(r) \cdot s_{S \times D^2}(u,r,\theta) + \sin^2 \Theta(r) \cdot
   s_{S \times S^1_r}(u,\theta) \geq s
\]
as long as $r \geq \mu + \epsilon$ (this is the region where 
$\lambda = \kappa$, so in particular the left hand side of this inequality is 
the scalar curvature of the metric $q$). Here we use the fact that 
$s_{(S \times ( D^2 \setminus 0), q)}$ is bounded below by the same constant $s$ 
as the scalar curvature of $(S \times (D(W)\setminus 0), h)$ -  but this metric on
$S \times (D(W) \setminus 0)$ of course extends to a positive scalar curvature metric 
on 
$S \times D(W)$ and consequently $s > 0$.  
Using the inequality $|\lambda''| \leq 2C$, 
which holds for any profile with bending parameter $C$, 
we know that for small enough  $\epsilon$, the angle $\Theta(r)$ changes so 
little in the region $r \in [\mu, \mu + \epsilon]$  that the previous sum 
is larger than $s/2>0$ for  all $r \geq \mu$. Because we can 
additionally assume that $\mu + \epsilon \leq  \overline{\mu}$, we conclude 
(by Lemma \ref{jan}) that the hypersurface $\Sigma$ 
has positive scalar curvature
for $r \geq \mu$ (this uses $k(r) \leq 0$ for any profile and any $r$). 
But once $r\leq \mu$, we are in the safe
region where the last term 
\[
   2k(r) \sin \Theta(r) \cdot K_{S \times S^1_r}(u,\theta)
\]
dominates - with a margin of at least $1$ - the sum of the absolute values 
of the first two terms in the expression of the 
scalar curvature of $\Sigma$ in Lemma \ref{jan}. 
Hence (and with  a continuity argument for $r=0$) with this profile $\lambda$,
the induced metric on $\Sigma$ is of positive scalar curvature. 

The projection 
\[
   \R^2 \to \R^2 \, , ~~ (r,t) \mapsto (r,0)
\]
induces a diffeomorphism 
\[
   \Sigma \approx S \times D^2 
\]
and we finally get an induced metric on $S \times D^2$ which is of 
positive scalar curvature and coincides near the boundary 
$S \times \partial D^2$      
with the metric on $S \times S(W)/H$ induced by $g$. 

In order to complete the proof of Theorem \ref{heikel}, we  choose a
trivialization of the $S^1$-principal bundle 
\[
    S^1 \hookrightarrow S^1 \times_H (S^{n-3} \times S(W)) \to 
      S^{n-3} \times S(W)/H \, .
\]
Now we pick an $S^1$-principal connection $\omega$ on the total space of the 
trivial $S^1$-principal bundle 
\[
   S^1 \hookrightarrow  S^1 \times (S^{n-3} \times D^2) \to S^{n-3} \times D^2
\]
which - after applying the above trivialization - 
coincides near the boundary $S^1 \times (S^{n-3} \times S^1)$ 
with the $S^1$-connection on the total space of
\[
    S^1 \hookrightarrow S^1 \times_H (S^{n-3} \times (D(W) \setminus 0)) \to S^{n-3} \times (D(W) \setminus 0)/H 
\]
which is induced by viewing $g$ as a Riemannian submersion metric on this fibre bundle.

For $\epsilon > 0$ we now consider the associated Riemannian submersion metric on 
\[
   S^1 \hookrightarrow S^1 \times (S^{n-3} \times D^2)
 \to S^{n-3} \times D^2
\]
with fibres $S^1$ of constant length $\epsilon$, horizontal 
subspaces induced by $\omega$  and the smooth positive scalar 
curvature metric on $S^{n-3} \times D^2$ constructed before. It follows from O'Neill that for small enough 
$\epsilon$, this metric has positive scalar curvature. Without loss of generality (possibly
after shrinking the orbits in $M \setminus N$) we can assume that the metric on $M \setminus N$ is of the same scale
(without violating the positive scalar curvature property). 
Hence the two metrics on $M \setminus N$ and on $S^1 \times (S^{n-3} \times D^2)$ can be combined such as to  
define an $S^1$-invariant metric of positive scalar curvature on $M'$. 
This completes the  proof of Theorem \ref{heikel}.

\section{Fixed point free $S^1$-manifolds} \label{anwend}

The main purpose of this section is to prove Theorem B from 
the introduction (see Theorem \ref{appl1} below). This proof 
is summarized as follows. A theorem of Ossa \cite{Ossa} states that 
any oriented fixed point free $S^1$-manifold $M$ satisfying condition $\C$ is 
the boundary of an oriented $S^1$-manifold $W$ (possibly with 
fixed points). However, our Theorem A cannot be applied directly because 
$W$ may contain singular strata of codimension $2$ 
that are disjoint from $M$. Different ideas are needed to 
handle this problem. At first, we remove tubular neighbourhoods 
of the components of the fixed point set  $W^{S^1}$ (which are disjoint from $M$ as 
$M^{S^1} = \emptyset$). This produces new boundary components of $W$ 
carrying invariant metrics of positive scalar curvature. This 
last statement follows 
from the O'Neill formula if the codimension of the corresponding  fixed component 
in $W$ is larger than $2$, and from the structure of the oriented bordism ring of 
free $S^1$-manifolds, which will be explained in the proof of Proposition \ref{free} below, 
if this codimension is equal to $2$ (we can 
assume that the action on $W$ is effective). Let $Z$ be the resulting fixed point free
bordism. Unfortunately, the remaining codimension-$2$ singular 
strata can be embedded in $Z$ in a complicated way. For example, their closures
may have nonempty intersections with each other (cf. Example \ref{typical}). The resulting 
problems can be circumvented by cutting out small equivariant tubes in $Z_{max}$ connecting $M$ with those
codimension-$2$ singular strata that are disjoint from $M$. In this way, 
we add new 
singular strata to $M$, but in any case Theorem A can be applied. 
Consequently,  the manifold $M'$ obtained from $M$ by adding these 
singular strata (and adding certain free $2$-handles, but we ignore 
this step for the moment) admits an $S^1$-invariant metric of positive scalar curvature. 
It turns out that $M$ can be recovered from $M'$ by performing 
surgery steps as explained in the previous Section \ref{genial}. In
particular,  thanks to Theorem \ref{heikel}, the original manifold $M$ 
admits an invariant metric of positive scalar curvature. The details 
of this argument are explained in the proof of Theorem \ref{seinsbesser}.

Let us start with the following observation concerning free $G$-manifolds. 

\begin{lem} \label{basic} Let $M$ be a closed manifold equipped with a 
free $G$-action. If the identity component of $G$ is abelian, then the following assertions are equivalent:
\begin{itemize} 
    \item[i.)] $M$ admits a $G$-invariant metric of positive scalar curvature. 
    \item[ii.)] $M/G$ admits a metric of positive scalar curvature. 
\end{itemize}
\end{lem}

\begin{proof} The case of finite $G$ is immediate. The general case
is Theorem C in \cite{Bebe}. 
\end{proof} 

Note that this fact is not true for connected nonabelian $G$. An 
easy counterexample is given by $M = \SU(2) \times S^1$ with a bi-invariant 
Riemannian metric on $\SU(2)$ (which has positive scalar curvature)
and $\SU(2)$ acting freely on the first factor in $\SU(2) \times S^1$.

Before we state the next proposition, we remind the reader of the 
following basic fact (cf.~\cite{CF}). Let $M$ and $N$ be closed 
oriented manifolds equipped with free 
orientation preserving $G$-actions. 
Then the following assertions are equivalent: 
 \begin{itemize}
     \item[i.)] There is a compact oriented  $G$-bordism $W$ between $M$ 
           and $N$ such 
           that $G$ acts freely and orientation preserving on $W$. 
     \item[ii.)] Consider the orbit manifolds $M/G$ and $N/G$ together 
           with the maps $f_M : M/G \to BG$ and 
           $f_N : N/G \to BG$ classifying the respective $G$-principal 
           bundles. Then $f_M: M/G \to BG$ and $f_N: N/G \to BG$ define 
           the same bordism class in $\Omega^{SO}_*(BG)$. 
 \end{itemize}

The following proposition contains our first general existence 
result of invariant metrics of positive scalar curvature on $S^1$-manifolds.

\begin{prop} \label{free} Let $M$ be a closed oriented free $S^1$-manifold of 
dimension at least $6$ which is simply connected and does not admit 
a spin structure. Then $M$ carries an $S^1$-invariant metric of positive 
scalar curvature. 
\end{prop}

\begin{proof} We give two proofs of this fact.

The  long exact homotopy sequence of the $S^1$-fibration
\[   
    S^1 \hookrightarrow M \stackrel{\pi}{\to} M/S^1
\]
shows that $M/S^1$ is also simply connected. Furthermore,
\[  
    TM \cong \pi^* (T (M/S^1)) \oplus \underline{\R}
\]
with a trivial line bundle $\underline{\R}$, and hence $M/S^1$ does 
not admit a spin structure by an easy  
characteristic class calculation. By the Gromov-Lawson theorem 
stated in the introduction,  $M/S^1$ admits a metric of positive scalar curvature and 
by Lemma \ref{basic}, the manifold $M$ admits an $S^1$-invariant 
metric of positive scalar curvature. 

The second proof is independent of \cite{Bebe} and a little longer, 
but can later be generalized to a wider class of actions.
 
By a standard use of the Atiyah-Hirzebruch spectral 
sequence, 
\[  
   \Omega^{SO}_*( B S^1) \cong \Omega^{SO}_*[x_1,x_2, \ldots ]
\]
is a graded polynomial ring in indeterminates $x_i$ of degree $2i$
(recall that $B S^1 = \C P^{\infty}$). The 
variable $x_i$ can be assumed to correspond to the free $S^1$-bordism 
class represented by the sphere 
\[
    S^{2i+1} \subset \C^{i+1}
\]
equipped with the standard free $S^1$-action. Obviously, these 
manifolds carry $S^1$-invariant metrics of positive 
scalar curvature. Together with the fact that each  
element in  $\Omega^{SO}_*$ can be represented
by a manifold admitting a positive scalar curvature metric 
(see \cite{GL}), this 
implies that each element in  $\Omega^{SO}_*(B S^1)$
is represented by a free $S^1$-manifold carrying an  invariant metric 
of positive scalar curvature. 

We conclude that the given manifold $M$ is bordant to an $S^1$-manifold admitting 
an invariant metric of positive scalar curvature and moreover the 
bordism $W$ can be assumed to be an oriented  free $S^1$-manifold. 
We need to show that $W$ can be improved in such a way that the inclusion 
$M \hookrightarrow  W$ is a $2$-equivalence. 

Because $\C P^{\infty}$ is 
simply connected, we can kill the fundamental group in $W / S^1$ 
by surgeries over (i.e. with reference maps to) $\C P^{\infty}$. 

Comparing  the long exact homotopy sequences induced by the commutative
diagram of fibrations
\[
   \begin{CD} 
      S^1 @>>> M @>>> M /S^1 \\
      @V = VV   @VVV @VVV    \\
      S^1 @>>> W @>>> W/S^1 
   \end{CD}
\]
we see that the new bordism $W$ is simply connected. 

Now let $c \in \pi_2(W/S^1)$ represent an element in the cokernel 
of the map 
\[
    \pi_2(M/S^1) \to \pi_2(W/S^1) \, . 
\]
We can represent $c$ by an embedded $2$-sphere $S^2 \subset W/S^1$ (recall 
that $\dim W/ S^1 \geq 5$). Let $\lambda$ be the image of $c$ under the map 
\[
   \pi_2(W/S^1) \to \pi_2(\C P^{\infty}) \cong \Z  
\]
which is induced by the reference map $W/S^1 \to \C P^{\infty}$. 
Before we can kill $c$ by surgery over $\C P^{\infty}$ we must make sure that 
$\lambda = 0$ (this would be automatic if we replaced $S^1$ by 
a finite group $G$ because then $\pi_2(BG) = 0$) 
and that the  normal bundle of $S^2 \subset W/S^1$
is trivial. 

In order to achieve these requirements we consider the commutative diagram  
\[
   \begin{CD}
     \pi_2(M/S^1) @>>> \pi_1(S^1) @>>> \pi_1(M) = \{ 1 \}  \\
         @VVV            @V =  VV \\
     \pi_2(W/S^1) @>>> \pi_1(S^1) \\
         @VVV            @V = VV \\
     \pi_2(\C P^{\infty}) @> \cong >> \pi_1(S^1)
   \end{CD}
\]
which is induced by the composition 
\[
    M/S^1 \hookrightarrow W/S^1 \to \C P^{\infty}
\] 
and whose horizontal maps are connecting homomorphisms in the 
long exact homotopy sequences of the respective $S^1$-fibrations. 
Because the first horizontal map is  surjective, we
find an element 
\[
    y \in \pi_2(M/S^1)
\]
which goes to $\lambda$ under the map 
\[
    \pi_2(M/S^1) \to \pi_1(S^1) \, . 
\]
After replacing $c$ by $c - y$ we can therefore assume that $\lambda = 0$.
If the second Stiefel-Whitney class of $W/S^1$ evaluated on (the new) $c$
is nontrivial,  we pick an element $x \in \pi_2(M)$ on 
which the second Stiefel Whitney class of $M$ evaluates
nontrivially (a spherical class in $H_2(M;\Z)$ 
with this property exists, because $M$ is not spin and simply connected). Now  
we replace $c$ by $c + \eta(x)$ where 
\[
   \eta : \pi_2(M) \to \pi_2(M/S^1) \to \pi_2(W/S^1)
\]
is the obvious map. This will preserve the property 
that $\lambda = 0$ because the composition 
\[
    \pi_2(M) \to \pi_2(M/S^1) \to \pi_1(S^1) 
\]
is zero. 

This shows that we can indeed kill $c$ by surgery 
over $\C P^{\infty}$. Because $\pi_2(W/S^1)$ is finitely
generated, we can therefore (after finitely many surgery steps) 
assume that the inclusion 
\[
   M/S^1 \hookrightarrow W/S^1
\]
is a $2$-equivalence and the same is then true for the inclusion 
$M \hookrightarrow W$. An application of Theorem \ref{possym} 
finishes the proof of Proposition \ref{free}. 
\end{proof}

Before we generalize the last proposition to fixed point free $S^1$-manifolds,
we show that assumption iii.) in Theorem \ref{possym} can be avoided
in the case of fixed point free $S^1$-actions whose union of 
maximal orbits is simply connected and not spin. Here the 
surgery procedure explained in Section \ref{genial} will be used.  

\begin{thm} \label{seinsbesser} Let $Z$ be a compact connected oriented fixed point 
free $S^1$-bordism between the closed $S^1$-manifolds $X$ and $Y$. Assume 
that $Z$ satisfies condition $\C$ and that the following hold:
\begin{itemize}
       \item[i.)] The cohomogeneity of $Z$ is at least $6$,  
       \item[ii.)] the union of maximal orbits $Y_{max}$ is simply connected 
                and does not admit a spin structure. 
\end{itemize}
Then, if $X$ admits an $S^1$-invariant metric of positive scalar
curvature which is normally symmetric in codimension $2$, 
the same is true for $Y$. 
\end{thm}

\begin{proof} Let $n = \dim X$ (i.e. $\dim Z = n+1$). 
By Lemma \ref{substantial1} we may assume that the given metric on $X$ is scaled. Now let  
\[
   F \subset Z 
\]
be a codimension-$2$ singular stratum in $Z$ (i.e. $\dim F = n-1$) which 
has empty intersection with $Y$ and is therefore problematic 
in view of Theorem \ref{possym}. By assumption, the 
isotropy group  $H \subset S^1$ of  $F$ is finite. Let $\Omega \subset  F$ 
be an orbit. It follows from the slice theorem that $\Omega$ 
has an $S^1$-invariant closed tubular neighbourhood $N$ in $Z$ which is  
$S^1$-diffeomorphic to 
\[
    S^1 \times_H (D^{n-2} \times D(W)) 
\]
where $W$ is a one dimensional unitary $H$-representation (because 
the given action on $Z$ satisfies condition $\C$) and 
$S^1$ acts only on the $S^1$-factor. We can assume that the $S^1$-action 
on $Z$ is effective and hence $H$ acts effectively on $W$. 
The idea is to alter $Z$ by cutting out an equivariant tube in $Z_{max}$ 
which connects  $N$ and  $Y$. 



We write
\[  
  \partial N = S^1 \times_H \big( ( D^{n-2} \times S(W)) \cup (S^{n-3} \times 
D(W) \big)  
\]
and use the  $H$-invariant subset  
\[
 T : =  \{ e^{\frac{2\pi i \omega}{|H|}} ~|~ \omega \in  \bigcup_{k=0}^{|H|-1} [k, k+\frac{1}{2 |H|}] \} \subset S(W) = S^1 
\]
to define the  $S^1$-invariant submanifold  
\[
    B := S^1 \times_H \big( D_{[0,1/2]}^{n-2} \times T) \subset \partial N  
\]
(the subscript at $D^{n-2}$ indicates restriction of the radial coordinate). 
The $S^1$-principal bundle 
\[
    S^1 \hookrightarrow B \to D_{[0,1/2]}^{n-2} \times T/H \approx D^{n-1} 
\]
is trivial. Hence, by the connectivity of $Z_{max}/S^1$, there exists
an orientation preserving 
$S^1$-equivariant embedding 
\[
   \Psi : (S^1 \times D^{n-1})  \times [0,1] \to Z_{max}
\]
(with $S^1$-acting freely on the $S^1$-factor) 
which restricts to an $S^1$-equivariant diffeomorphism 
\[
     S^1 \times D^{n-1} \times \{0\} \approx B \subset \partial N  
\]
and to an embedding 
\[
     S^1 \times D^{n-1} \times \{1\} \hookrightarrow Y 
\]
and satisfies
\[
    \Psi\big( S^1 \times D^{n-1} \times (0,1)\big) \subset Z \setminus (Y \cup N) \, .
\]
We now consider the $S^1$-bordism 
\[
  Z' :=  Z \setminus \big( N \cup \im (\Psi)\big) \, . 
\]
In this bordism, the manifold $Y$ is replaced by another manifold $Y'$ 
which contains a new codimension-$2$ singular stratum 
\[  
   \Sigma :=  S^1 \times_H (S^{n-3} \times 0) \subset S^1 \times_H 
              (S^{n-3} \times D(W)) \subset \partial N  \, . 
\]
We claim that $Y$ can be recovered from $Y'$ by resolving $\Sigma$.  
The argument goes as follows: The construction of $Y'$ yields an embedding 
\[
   \phi' :   S^1 \times_H (S^{n-3} \times D(W)) \hookrightarrow Y'
\]
of a tubular neighbourhood of $\Sigma \subset Y'$ and the manifold 
$Y' \setminus \im(\phi')$ can be written as 
\[
   \big(  Y \setminus \Psi (S^1 \times D^{n-1} \times \{1\}) \big) 
   \cup_{S^1 \times S^{n-2} \times \{1\}} \big( S^1 \times \partial D^{n-1}  \times [0,1]) \cup_{\partial B}  A  
\]
where    
\[
 A :=  S^1 \times_H \big( (D^{n-2} \times S(W)) \setminus (D^{n-2}_{[0, 1/2]} \times T)  \big)  \subset \partial N \, . 
\]
The $S^1$-action on $A$ is free and the quotient space 
$A/S^1$ is diffeomorphic to $(D^{n-2} \times S^1) \setminus D$
where $D= B/S^1$ is a submanifold of $D^{n-2} \times S^1$ diffeomorphic to $D^{n-1}$.  
Because $n \geq 6$ by assumption, all principal $S^1$-bundles 
over $A/S^1$ are isomorphic and hence there is an $S^1$-equivariant 
diffeomorphism. 
\[
   A \approx S^1 \times \big( (D^{n-2} \times S^1) \setminus D \big) \, . 
\]
We conclude that there is an $S^1$-equivariant diffeomorphism 
\[
   Y' \setminus \im (\phi') \approx  Y \setminus \im (\phi) 
\]
where 
\[
   \phi : S^1 \times (S^{n-3} \times D^2) \hookrightarrow Y_{max} 
\]
is an $S^1$-equivariant orientation preserving embedding whose
image is contained in the $S^1$-equivariant coordinate
chart 
\[
     \Psi(S^1 \times D^{n-1} \times \{1\}) \subset Y \, . 
\]
(Note the standard decomposition $S^{n-1} =  (D^{n-2} \times S^1)\cup (S^{n-3} \times D^2)$.)
It follows that we can write $Y$ as 
\[
 \big(  Y' \setminus {\rm im} (\phi') \big) \cup \big( S^1 \times (S^{n-3} \times
          D^2) \big) 
\]
and this proves that $Y$ can be recovered from $Y'$ by a resolution of 
$\Sigma$.    

In particular (using Theorem \ref{heikel}), 
if we can show that $Y'$ admits 
a scaled  $S^1$-invariant metric of positive scalar curvature which is 
normally symmetric in codimension $2$, the same holds for $Y$.  

Because the embedded $S^1$-manifold 
\[
   \phi( S^1 \times (S^{n-3} \times 0) ) \subset Y_{max}
\]
is contained  in an $S^1$-equivariant coordinate chart, it 
can be assumed to be disjoint from some embedded 
$2$-sphere in $Y_{max}$ with 
nontrivial normal bundle (such a $2$-sphere exists because $Y_{max}$ does 
not admit a spin structure). This implies that $Y'_{max}$ does 
not admit a spin structure, either.  We would like $Y'_{max}$ 
to be simply connected, too. However, this need not be the
case due to the existence of a non-nullhomotopic linking sphere 
\[
     S^1 \subset Y'_{max}/S^1
\]
of $\Sigma/S^1 = S^{n-3} \subset Y'_{max}$. But this problem can 
be solved as follows: Before we perfom the cutting-out procedure on $Z$, 
we attach a free $S^1$-equivariant $2$-handle $S^1 \times (D^2 \times D^{n-2})$ to 
\[
    S^1 \times (S^{n-3} \times D^2) \subset Y_{max} \subset \partial Z 
\]
(here we suppress the identification $\phi$) along 
\[
   S^1 \times (S^1 \times D^{n-2}) \subset S^1 \times (S^{n-3} \times D^2)  
\]
where $\{1\} \times (S^1 \times D^{n-2})$ is identified with 
a small tubular neighbourhood of 
$\{1\} \times \{ p \} \times S^1_{1/2} \subset \{1\} \times S^{n-3} 
\times D^2$. Here, 
$p \in S^{n-3}$ is an arbitrary point and $S^1_{1/2} \subset 
D^2$ is the circle of radius  $1/2$. The space which 
is obtained from $Z$ by attaching this free $2$-handle is denoted
by $\widetilde{Z}$. By construction, we can attach a 
further free $S^1$-equivariant $3$-handle 
to $\widetilde{Z}$ which may be canceled against the previously 
attached $2$-handle. Furthermore (by a backward use of 
the Seifert-van Kampen theorem)  
\[
  \pi_1 ( \widetilde{Y} \setminus \big( S^1 \times (S^{n-3} \times D^2_{[0,1/4]}) 
  \big) = \{ 1 \}
\]
where $\widetilde{Y}$ is the space obtained from $Y$ by performing 
the surgery associated to the additional free $2$-handle (this 
can be assumed not to affect the subset $S^1 \times (S^{n-3} \times 
D^2_{[0, 1/4]}) \subset Y$). The old cutting 
out process on $Z$ can also be performed on the new bordism $\widetilde{Z}$
because it can be assumed only to affect  the part
\[
    S^1 \times (S^{n-3} \times D^2_{[0,1/4]}) \subset \widetilde{Y} \, .  
\]
The same procedure (i.e.~attaching a free 
dummy $2$-handle and  cutting out 
a certain part of the bordism)  is now applied
to all other singular strata of codimension $2$ in $Z$. In this way, we end 
up with an $S^1$-bordism $Z'$ in which all singular strata of codimension $2$ have 
nonempty intersection with $Y'$ and $Y'_{max}$ is simply connected and does 
not admit a spin structure.

We now attach equivariant handles to $Z'_{max}$  as in the proof of Proposition 
\ref{free} to make sure that the inclusion 
\[
     Y'_{max} \hookrightarrow Z'_{max}
\]
is $2$-connected. 
Here we note that $Z'_{max}$ has
finitely generated fundamental group and homology groups 
so that in any case, only finitely many surgery steps 
on $Z'_{max}/S^1$ are needed.

Theorem \ref{possym} together with its refinement 
formulated in Lemma \ref{substantial2} implies that $Y'$ admits 
a scaled $S^1$-invariant metric of positive scalar curvature
which is normally symmetric in codimension $2$.  
By Theorem \ref{heikel}, this also holds for 
the manifold obtained from $Y'$ by resolving the singularities
created by the cutting out processes. 
But then, as explained before, the attachment of the dummy  
free $2$-handles (which we can assume to be disjoint from each other)
can be neutralized by attaching free $3$-handles. 
These handles have codimension $\coh(Y,S^1)-2$ which is at least $3$ by assumption. Therefore, 
using the equivariant surgery principle, Theorem \ref{equivgl},  
this step also preserves the $S^1$-invariant scalar curvature metric and 
the resulting space - which can be identified with $Y$ - indeed 
carries an $S^1$-invariant metric of positive scalar curvature. 
\end{proof}

The following is Theorem B from the introduction.

\begin{thm} \label{appl1} Let $M$ be a closed
fixed point free $S^1$-manifold satisfying condition $\C$ and 
of cohomogeneity at least $5$. 
If $M_{max}$ is simply connected and does 
not admit a spin structure, then 
$M$ admits an $S^1$-invariant metric of positive scalar curvature. 
\end{thm}

\begin{proof} Without loss of generality, the given $S^1$-action on $M$ 
is effective. Because $M$ satisfies condition $\C$, the singular strata in $M$ are 
of codimension at least $2$ and hence $M$ is simply connected by a general 
position argument. In particular, 
it is an orientable $S^1$-manifold. 
By \cite{Ossa}, Satz 1,   the manifold $M$ is the boundary 
of an oriented (connected) $S^1$-manifold $W$ satisfying condition $\C$. Let 
\[
   F_1, \ldots, F_k
\]
be the components of $W^{S^1}$. By assumption, these are disjoint from $M$.
We cut out pairwise disjoint $S^1$-invariant tubular 
neighbourhoods $N_i$ of $F_i$ in $W$. This yields an oriented fixed point 
free bordism satisfying condition $\C$ from $M$ to another 
$S^1$-manifold with components $\partial N_i$, $1 \leq i \leq k$. 
Each $N_i$ is the total 
space of a unitary $S^1$-equivariant fibre bundle
\[
   V_i \hookrightarrow N_i \to F_i
\]
with a unitary $S^1$ representation $V_i$. If the codimension of $F_i$
is larger than $2$, then $\partial N_i$ carries an $S^1$-invariant metric
 of positive scalar curvature which 
is normally symmetric in codimension $2$ by the O'Neill 
formula in combination with $\dim S(V_i) \geq 2$  and 
Example \ref{oversym}. 
If ${\rm codim}\, F_i =2$, then, because the action on $W$ is effective,
the $S^1$-action on $V_i$ is effective 
and therefore the induced action on $\partial N_i$ is 
free. As explained in 
the second proof of Proposition \ref{free}, 
$\partial N_i$ is then freely and oriented bordant to a free $S^1$-manifold admitting 
an invariant metric of positive scalar curvature. 
Now Theorem \ref{appl1} follows from Theorem \ref{seinsbesser}. 
\end{proof}

One might ask whether Theorem \ref{appl1} can be proven without 
the somewhat involved discussion of codimension-$2$ singular 
strata in Section \ref{genial} and Theorem \ref{seinsbesser}
under the assumption that  
$M$ does not contain such strata. But a closer look at \cite{Ossa}
reveals that in general (depending on the dimensions of 
the isotypical summands of the normal representations  around the 
singular strata in $M$) the zero bordism $W$ does contain 
codimension-$2$ singular strata with finite isotropies, even if $M$ does not.

Unfortunately, we do not have such a general existence result for $S^1$-manifolds with fixed points. 
One can check that the oriented $S^1$-bordism 
ring (always restricting to actions satisfying condition $\C$) 
is generated by $S^1$-manifolds admitting $S^1$-invariant metrics 
of positive scalar curvature. This follows from an 
inspection of the generators constructed in  \cite{KY}. However, if 
we express a given $S^1$-manifold without codimension-$2$ singular 
strata in terms of these generators, it might happen that 
generators with codimension-$2$ singular strata do appear and this 
leads to $S^1$-handle decompositions of the given $S^1$-bordism 
containing handles of codimension $0$ or $2$, cf. Proposition \ref{illustrate}. 
However, we do not know if a surgery principle as explained in 
Section \ref{genial} exists for $S^1$-manifolds with fixed points 
(note that invariant metrics on such manifolds can never be 
scaled). 

In some special situations, one can construct the 
necessary bordisms by hand. For example, we have the 
following result for semifree $S^1$-manifolds with isolated 
fixed points.

\begin{thm} Let $M$ be a closed simply connected non-spin manifold of 
even dimension at least $6$ and equipped with a semifree $S^1$-action 
(i.e. the action has either free or fixed orbits) with 
only isolated fixed points. Then $M$ admits an $S^1$-invariant metric 
of positive scalar curvature.  
\end{thm} 

\begin{proof} Let $2n$ be the dimension of $M$. 
After removing small invariant discs around the 
fixed points and dividing out the free $S^1$-action, we
get a zero bordism over $B S^1$ of a disjoint union of copies of 
$\pm \C P^{n-1}$ where the reference maps to $BS^1$ classify the 
tautological line bundle over $\C P^{n-1}$. This 
classifying map $\C P^{n-1} \to BS^1$
generates a $\Z$-summand 
in $\Omega^{SO}_{2n-2}( B S^1)$ and  therefore we get as many $-$-signs
as we get $+$-signs. By pairwise connecting a positively oriented 
fixed point with a negatively oriented one by thin 
tubes we obtain an oriented $S^1$-bordism $W$
from a free  $S^1$-manifold $N$ to the given manifold.
Using the structure of $\Omega^{SO}_*(B S^1)$ (cf. the second 
proof of Proposition \ref{free})
we can  assume (possibly after adding a free oriented $S^1$-bordism to $N$)  
that  $N$ 
has an invariant metric
of positive scalar curvature.
Furthermore, 
by an argument similar to the second proof of  Proposition \ref{free}, 
the inclusion $M_{max} \hookrightarrow W_{max}$ can be assumed to be a 
$2$-equivalence. 
Theorem \ref{possym} now implies that the manifold $M$ admits 
an $S^1$-invariant metric of positive scalar curvature because 
assumption iii.) obviously holds. 
\end{proof}

\end{document}